\documentclass[10pt]{amsart}

\usepackage[leqno]{amsmath}
\usepackage{amssymb}
\usepackage{amsthm}
\usepackage{amscd}

\swapnumbers\newtheorem{thm}{Theorem}[section]\newtheorem{cor}[thm]{Corollary}\newtheorem{prop}[thm]{Proposition}\newtheorem{lem}[thm]{Lemma}\theoremstyle{definition}\newtheorem*{defn}{Definition}\newtheorem*{rem}{Remark}\newtheorem*{rems}{Remarks}
\newtheorem*{thmA}{Theorem A}\newtheorem*{thmB}{Theorem B}\newtheorem*{thmC}{Theorem C}\newtheorem*{thmD}{Theorem D}
\newtheorem*{thmE}{Theorem E}\newtheorem*{thmF}{Theorem F}\newtheorem*{thmG}{Theorem G}\newtheorem*{thmH}{Theorem H}
\newtheorem*{thmJ}{Theorem J}\newtheorem*{thmK}{Theorem K}
\newtheorem*{corI}{Corollary I}
\newtheorem*{corL}{Corollary L}\newtheorem*{corM}{Corollary M}
\newtheorem*{corN}{Corollary N}\newtheorem*{corO}{Corollary O}\newtheorem*{corP}{Corollary P}\newtheorem*{corQ}{Corollary Q}

\newcommand{\set}[1]{\left\{#1\right\}}
\newcommand{\abs}[1]{\left\vert#1\right\vert}
\newcommand{\Irr}{\mathrm{Irr}}

\newcommand{\core}{\mathrm{Core}}
\newcommand{\hm}{\mathrm{Hom}}

\newcommand{\res}{\mathrm{res}}
\newcommand{\Aut}{\mathrm{Aut}}

\begin{document}

\title{Primitive characters and permutation characters of solvable groups}
\author{Tom Wilde}

\begin{abstract}
If $\chi$ is a primitive character of the solvable group $G,$ and if $\chi(1)$ is odd, then we associate to $\chi$ in a unique way, a conjugacy class of subgroups $U\subseteq G$ which satisfy $\bar\chi\chi=(1_U)^G.$ Furthermore, if $G$ has odd order and $\bar\chi\chi=(1_U)^G$ for some subgroup $U\subseteq G,$ then conditions on $U$ exist which are sufficient for $\chi$ to be primitive. We investigate this, giving some applications to properties of primitive characters of solvable groups.
\end{abstract}

\maketitle

\section{Introduction}
Let $G$ be a finite group. Then $\chi\in\Irr(G)$ is \emph{quasi-primitive} if its restriction to any normal subgroup $N\vartriangleleft G$ is a multiple of an irreducible character of $N,$ and $\chi$ is \emph{primitive} if it cannot be obtained by inducing a character of a proper subgroup of $G.$ In any group, Clifford theory shows that primitive characters are quasi-primitive, but for solvable groups the converse is also true by Berger's well known theorem (\cite[Theorem 11.33]{Isaacs1976}). We are always concerned with solvable or $p$-solvable groups, but we usually refer to quasi-primitivity in the statements of results, to indicate where quasi-primitivity is the property actually used in a proof.\\

This paper arose from considering an interesting theorem, due to Gabriel Navarro, on the zeros of primitive characters of $p$-solvable groups (see Theorem A in Section 2, below), and afterwards reading the paper \cite{FergusonTurull1985} of Ferguson and Turull, which presents a theory of factorization of quasi-primitive characters. These sources lead us to investigate permutation characters associated with quasi-primitive characters of suitably solvable groups. If $G$ is solvable, and $\chi\in\Irr(G)$ is quasi-primitive and has odd degree, then we show that there always exists a subgroup $U\subseteq G$ such that \begin{equation*}\bar\chi\chi=(1_U)^G,\tag{\dag}\end{equation*} the permutation character of $G$ acting on the cosets of $U.$ Moreover, $U$ can be chosen in a prescribed manner. See Theorem C in Section 2. There exist irreducible characters of solvable groups which are not quasi-primitive, and yet which have an expression of the form (\dag). However, at least when $G$ is odd, conditions on the subgroup $U$ in (\dag) may be found which are necessary and sufficient for $\chi$ to be quasi-primitive. These conditions are the subject of Theorems E and F. The existence of such conditions implies that in a group $G$ of odd order, if $\chi,\psi\in\Irr(G)$ and $\bar\chi\chi=\bar\psi\psi,$ (i.e. if $\abs{\chi(g)}^2=\abs{\psi(g)}^2$ for all $g\in G$), then $\chi$ is primitive if and only if $\psi$ is primitive. (See the remarks after Theorem F in Section 2.)

Using our results, we give a nice proof of a theorem of Isaacs on the restriction of certain quasi-primitive characters to normal subgroups, in the case $\abs{G}$ odd (see Corollary I), and we show that an analogous result holds in some cases for non-normal subgroups (Theorem J). We also give a proof, in the case that $G$ has odd order, of a result of Ferguson and Isaacs (their result is Theorem A in \cite{FergusonIsaacs1989}), showing that multiples of primitive characters also cannot be induced from proper subgroups (see Theorem K)\footnote{The prior result of Ferguson and Isaacs, which includes Theorem K, was kindly pointed out to me by I. M. Isaacs. I regret the lack of attribution in the original version of this paper.}. In Corollaries L - Q, we give various properties of characters of solvable groups which are direct consequences of the existence of an expression (\dag).

The paper is organized as follows. Section 2 contains statements of our results. In Section 3, we study \emph{strongly irreducible} characters (see that section for the definition), and prove Theorems A and B. Our proofs of these results are elementary, apart from a standard application of projective representations. Section 4 is concerned with (mostly known) properties of solvable groups, and it is convenient to prove Theorems G and H there. Then in Section 5, we prove Theorems C and D. Section 6 contains the proofs of Theorems E, F, J and K. The short proofs of Corollaries I and L - Q are given immediately after their statements in Section 2.

We have attempted in this paper to give a fairly self-contained development, indicating where results are known. In particular, part of the work can be circumvented by appealing to a deep result of Isaacs (\cite[Theorem 9.1]{Isaacs1973}) on ramified sections of odd order. We also need a consequence of this well known result for our proof of Theorem K. See Theorem \ref{bigthm} and Corollary \ref{bigcor} below. There is also a degree of overlap with the paper \cite{FergusonTurull1985} of Ferguson and Turull referred to above. 

Groups considered in this paper are finite. We have tried to state our results with only the necessary solvability conditions. These always include the condition that $G$ is $p$-solvable for each prime $p$ dividing $\chi(1),$ i.e. $G$ is $\pi(\chi(1))$-solvable, where $\pi(n)$ denotes the set of prime divisors of $n.$ The reader who prefers to, may assume that all groups considered are solvable and ignore the various solvability conditions in our definitions and results. We usually need $\chi(1),$ or even $\abs{G:Z(\chi)}$ to be odd, where as usual $Z(\chi)$ is the pre-image in $G$ of the centre of $G/\ker\chi.$ 

\section{Statement of results}
The following is the theorem of Navarro mentioned in the introduction.

\begin{thmA}[\cite{Navarro1999}, Theorem A]
Let $p$ be a prime and suppose $G$ is $p$-solvable. Suppose $\chi\in\Irr(G)$ is quasi-primitive and has $\chi(1)$ a power
of $p.$ Let $x\in G,$ and let $x_p$ denote the $p$-part of $x.$ Then $\chi(x)\neq 0$ if and only if $\chi(x_p)\neq 0.$
\end{thmA}

Our proof also yields:

\begin{thmB}
Suppose $\chi\in\Irr(G)$ is quasi-primitive, and assume $G$ is $\pi(\chi(1))$-solvable. Then $\abs{\chi(x)}^2$ is either zero, or is a
positive integer dividing $\chi(1)^2.$
\end{thmB}

Our main results are Theorems C - K below. To begin, we need the following (non-standard) definition.

\begin{defn}
Let $U\subseteq G,$ be a subgroup of the group $G.$ Then $U$ is a \emph{complete intersection} if the following hold.
\begin{enumerate}
\item $G$ is $p$-solvable for each prime $p$ dividing $\abs{G:U}.$
\item There exist maximal subgroups $M_i\subset G$ such that $$U=\bigcap_iM_i\text{ and }\abs{G:U}=\prod_i\abs{G:M_i}.$$
\end{enumerate}

As mentioned in the introduction, the reader who prefers can assume that all groups considered are solvable and ignore condition (1) and analogous conditions below.

If the $M_i$ forming the intersection in condition (2) are to be specified, then we say that $U$ is the complete intersection \emph{of the} $M_i.$ We regard
$G$ itself as the complete intersection of the empty set of maximal subgroups. If $U\subseteq H\subseteq G$ then to avoid ambiguity we
may say $U$ is a complete intersection \emph{in $G$} or \emph{in $H$}, to indicate that the $M_i$ in condition (2) are
maximal subgroups of $G,$ or respectively of $H.$
\end{defn}

\begin{thmC} Suppose $\chi\in\Irr(G)$ is quasi-primitive and $\chi(1)$ is odd, and assume $G$ is $\pi(\chi(1))$-solvable.
Let $\mathfrak M(\chi)$ be the set of subgroups $H$ of $G$ which satisfy that $H$ is a complete intersection, and that $\chi\vert_H$ has no
zeros. Let $\mathfrak M_0(\chi)$ be the set of maximal elements of $\mathfrak M(\chi),$ with respect to its partial ordering by inclusion. Then $\mathfrak M_0(\chi)$ consists of exactly one conjugacy class of subgroups of $G,$ and for $U\in\mathfrak M_0(\chi),$ we have $\bar\chi\chi=(1_U)^G.$
\end{thmC}

\begin{proof}
See Section 5.
\end{proof}

Theorem C does not hold in general when $\chi(1)$ is even. See Remark (1) at the end of this section. The next result expresses a relationship
between Theorem C and the factorizations of quasi-primitive characters studied by Ferguson and Turull in \cite{FergusonTurull1985}. Here and often in the sequel, we use the concept of an \emph{extension} of a group $G,$ which means a group $\tilde G$ and a homomorphism $\pi:\tilde G\twoheadrightarrow G$ onto $G.$ The extension is \emph{central} if $\ker\pi\subseteq Z(\tilde G).$ 

\begin{thmD}
Suppose $\chi\in\Irr(G)$ is quasi-primitive and assume $\chi(1)$ is odd and $G$ is $\pi(\chi(1))$-solvable. By Theorem C, let
$U\in\mathfrak M_0(\chi).$ Then there exist pairwise inconjugate maximal subgroups $M_i$ of $G,$ a central extension $\pi:\tilde
G\twoheadrightarrow G,$ and irreducible characters $\rho_i\in\Irr(\tilde G)$ such that
\begin{enumerate}
\item $U$ is the complete intersection of the $M_i.$
\item $\bar\rho_i\rho_i=(1_{\tilde M_i})^{\tilde G},$ where $\tilde M_i=\pi^{-1}(M_i),$ the pre-image of $M_i$ in $\tilde G.$ 
\item $\tilde\chi=\prod_i\rho_i,$ where $\tilde\chi=\chi\circ\pi$ is the lift of $\chi$ to $\tilde G.$
\end{enumerate}
\end{thmD}
\begin{proof}
See Section 5.
\end{proof}

Our next results come from asking what properties characterize the subgroups $U$ in the statement of Theorem C. If $\chi\in\Irr(G)$ where $G$ is solvable and $\bar\chi\chi=(1_U)^G$ for a subgroup of $G,$ then we cannot conclude that $\chi$ is quasi-primitive without further conditions on $U,$ even if $U$ is a complete intersection as in Theorem C. For example, if $\chi_i,$ $i=1,2$ are quasi-primitive and satisfy $\bar\chi_i\chi_i=(1_{U_i})^G$ for subgroups $U_1,U_2,$ and if $\chi=\chi_1\chi_2$ is irreducible, then it is easy to see that $\bar\chi\chi=(1_{U_1\cap U_2})^G$ and $U_1\cap U_2$ is a complete intersection. However, $\chi_1\chi_2$ need not be quasi-primitive.

To state Theorems E and F, we introduce some further notation. Suppose $M$ is a maximal subgroup of the $p$-solvable group $G,$ and suppose $p$
divides $\abs{G:M}.$ Then $G/\core_G(M)$ has a unique minimal normal subgroup, which we will always denote by $X_M.$ In fact, $X_M$ is an elementary abelian $p$-group, and we regard it as a simple right $\mathbb{F}_pG$-module via conjugation. See Lemma \ref{XMlem} for these facts. Also, for
any abelian group $X,$ acted on by $G,$ we let $X^*$ denote the dual module, i.e. $\Irr(X)$ with natural $G$-action $\lambda^g(x^g)=\lambda(x).$

\begin{defn}
Let $G$ be a finite group and let $U\subseteq G$ be a subgroup of $G.$ Then $U$ is a \emph{regular intersection} if there exist maximal
subgroups $M_i\subset G$ such that the following conditions hold:
\begin{enumerate}
\item $G$ is $p$-solvable for each prime $p$ dividing $\abs{G:U}.$
\item $U=\bigcap_iM_i.$
\item For any $i\neq j,$ $X_{M_i}$ is isomorphic neither to $X_{M_j},$ nor to $X_{M_j}^*$ as $G$-module.
\end{enumerate}
\end{defn}

A regular intersection is, in particular, a complete intersection (see Lemma \ref{rcilemma}(1) in Section 3). The following pair of
theorems shows that, at least in a group of odd order, we have a characterization of quasi-primitive characters, as exactly the
irreducible characters which satisfy $\bar\chi\chi=(1_U)^G$ for a regular intersection $U\subseteq G.$ Actually, only one direction (Theorem F) requires an oddness assumption.

\begin{thmE}
Suppose $\chi\in\Irr(G)$ satisfies $\bar\chi\chi=(1_U)^G,$ for a subgroup $U\subseteq G.$
If $U$ is a regular intersection, then $\chi$ is quasi-primitive.
\end{thmE}
\begin{proof}
See Section 6.
\end{proof}

Note that in view of part (1) of the definition of a regular intersection, the group $G$ in Theorem E is automatically $\pi(\chi(1))$-solvable.

\begin{thmF}
Suppose $\chi\in\Irr(G)$ is quasi-primitive, and assume $G$ is $\pi(\chi(1))$-solvable and $\abs{G:Z(\chi)}$ is odd (so that in
particular, $\chi(1)$ is odd and Theorem C applies). Let $U\in\mathfrak M_0(\chi).$ Then $U$ is a regular intersection.
\end{thmF}
\begin{proof}
See Section 6.
\end{proof}

\begin{rems}
1. The concept of a regular intersection is parallel to the odd order case of the \emph{admissible} sets of prime characters introduced by Ferguson and Turull in \cite{FergusonTurull1985}.

2. Theorem F does not exclude the possibility that $\chi$ has an ``extra representation'' $\bar\chi\chi=(1_V)^G$ for some subgroup $V\notin\mathfrak M_0(\chi).$ If so, then Theorem C shows that $V$ cannot be a complete intersection. We do not know if such an ``extra representation'' can actually occur, however.

3. As remarked in the Introduction, Theorems E and F imply that, if $\abs{G}$ is odd and $\chi,\psi\in\Irr(G)$ have $\bar\chi\chi=\bar\psi\psi,$ then $\chi$ is quasi-primitive if and only if $\psi$ is quasi-primitive. For if one of $\chi$ or $\psi$ is quasi-primitive, say $\chi,$ then by Theorem F, $U\in\mathfrak M_0(\chi)$ is a regular intersection and satisfies $\bar\chi\chi=\bar\psi\psi=(1_U)^G,$ so both $\chi$ and $\psi$ are quasi-primitive by Theorem E.
\end{rems}

In the situation of Theorem C, if $U\in\mathfrak M_0(\chi)$ is a regular intersection (as is the case when $\abs{G:Z(\chi)}$ is odd, by Theorem F), then $U$ also has a property dual to that of Theorem C. Namely, if $V\subseteq G$ is a complete intersection, then $\chi$ vanishes off the union of the conjugates of $V$ if and only if $V$ is conjugate to a subgroup of $U.$ Thus in this case, $\mathfrak M_0(\chi)$ is also characterized as the set of subgroups of $G$ which are complete intersections, have $\chi$ supported on the union of their conjugates, and are minimal subject to these
conditions. This is a consequence of Theorem G below.

\begin{thmG}
Suppose $\chi\in\Irr(G)$ and $\bar\chi\chi=(1_U)^G$ where $U$ is a regular intersection. (This holds by Theorem F when $\chi$ is primitive and $\abs{G:Z(\chi)}$ is odd.) Let $V\subseteq G$ be a complete intersection. Then the following are equivalent
\begin{enumerate}
\item $\chi$ vanishes off the union of the $G$-conjugates of $V.$
\item $U$ is conjugate to a subgroup of $V.$
\end{enumerate}
Furthermore, if $V$ is a maximal subgroup of $G,$ then either (1) and (2) hold, or $\chi\vert_V$ is irreducible.
\end{thmG}
\begin{proof}
See Section 4.
\end{proof}

Theorems F and G imply the following fact about restricting quasi-primitive characters to maximal subgroups: If $\abs{G}$ is odd, $V$ is a maximal subgroup of $G,$ and $\chi\in\Irr(G)$ is quasi-primitive, then either $\chi$ vanishes off the union of the $G$-conjugates of $V,$ or $\chi\vert_V$ is irreducible. 

In the rest of this section we give some applications of Theorems C - F. First, in Theorem H and Corollary I, we give a proof of a theorem of Isaacs
(\cite[Theorem B]{Isaacs1981}) in the case $\abs{G}$ odd.

\begin{thmH}
Let $p$ be any prime. Suppose $N\vartriangleleft G$ and $\abs{G:N}$ is a power of $p.$ Suppose $U\subseteq G$ is a
regular intersection such that $\abs{G:U}$ is also a power of $p.$ Then $U\cap N$ is a regular intersection in $N.$
\end{thmH}

\begin{proof}
See Section 4.
\end{proof}

\begin{corI}
Let $G$ be $p$-solvable and suppose $\abs{G}$ is odd. Let $N\vartriangleleft G$ with $\abs{G:N}$ a power of $p.$ Suppose
$\chi\in\Irr(G)$ is quasi-primitive and has $\chi(1)$ also a power of $p.$ Then $\chi\vert_N$ is irreducible and quasi-primitive, and
$\mathfrak M_0(\chi\vert_N)=\set{U\cap N\vert U\in\mathfrak M_0(\chi)}.$
\end{corI}

\begin{proof}
By Theorem C, $\bar\chi\chi=(1_U)^G$ where $U\in\mathfrak M_0(\chi),$ and by Theorem F, since $\abs{G}$ is odd, 
$U$ is a regular intersection in $G.$ If $NU$ is a proper subgroup of $G,$ then since $G/N$ is nilpotent, there is a normal subgroup $M\vartriangleleft G$ with $\abs{G:M}=p$ and $NU\subseteq M.$ Then since $\bar\chi\chi=(1_U)^G,$ it is clear that $\chi$ vanishes off $M,$ and since $M$ is normal of prime index, therefore $\chi$ is induced from $M,$ a contradiction. We conclude that $NU=G.$ Thus $[\chi,\chi]_N=[(1_U)^G,1]_N=1,$ so $\chi\vert_N$ is irreducible and $\bar\chi\chi\vert_N=(1_{U\cap N})^N.$ By Theorem H, $U\cap N$ is a regular intersection in $N,$ so by Theorem E, $\chi\vert_N$ is quasi-primitive. Since $U\cap N$ is a regular intersection, it is a complete intersection (we prove this easy fact at Lemma \ref{rcilemma}(1) in Section 4), and since $\bar\chi\chi\vert_N=(1_{U\cap N})^N,$ clearly $\chi\vert_{U\cap N}$ has no zeros, so $U\cap N\in\mathfrak M(\chi\vert_N).$ Then by Theorem C, there exists $W\in\mathfrak M_0(\chi\vert_N)$ with $U\cap N\subseteq W,$ and equality holds since $\chi(1)^2=\abs{G:W}$ by Theorem C. Hence $U\cap N\in\mathfrak M_0(\chi\vert_N)$ as required.
\end{proof}

In Corollary I, the subgroup $N$ is a normal subgroup of $G.$ The following curious result shows that the normality assumption can be dispensed with when $\chi(1)$ is not too large.

\begin{thmJ}
Suppose $\chi\in\Irr(G)$ and $\bar\chi\chi=(1_U)^G$ where $U$ is a regular intersection. Assume $\chi(1)=p^a,$ where $p$ is an odd prime and $a<p.$ Let $H$ be a subgroup of $G$ with $\abs{G:H}$ also a power of $p.$ If $\chi\vert_H$ is irreducible, then it is quasi-primitive.
\end{thmJ}

\begin{proof}
See Section 6.
\end{proof}

Our next result is a strengthened form of Berger's theorem (\cite[Theorem 11.33]{Isaacs1976}) for groups of odd order. Berger's theorem for groups of odd order is the case $e=1.$ After the original version of this paper was written, it was kindly pointed out to us by I. M. Isaacs that this result is essentially the odd order case of a theorem of Ferguson and Isaacs, Theorem A in \cite{FergusonIsaacs1989}.

\begin{thmK}[cf. Theorem A in \cite{FergusonIsaacs1989}]
Let $\chi\in\Irr(G)$ be quasi-primitive and suppose $G$ is $\pi(\chi(1))$-solvable and $\abs{G:Z(\chi)}$ is odd. Suppose $H\subseteq G$ and $\theta$ is a character of $H$ with $\theta^G=e\chi$ for some $e\in\mathbb N.$ Then $H=G.$ 
\end{thmK}

\begin{proof}
See Section 6.
\end{proof}

\begin{corL}
Let $N\vartriangleleft G$ where $G$ is solvable, and assume $N$ has odd order. Let $\chi\in\Irr(G)$ and suppose $\chi\vert_N=e\theta$ where
$\theta\in\Irr(N)$ is quasi-primitive. Suppose $\chi=\mu^G$ for a subgroup $H\subseteq G$ and $\mu\in\Irr(H).$ Then $H\supseteq N$ and
$\mu_N=f\theta.$
\end{corL}

\begin{proof}
By Mackey's theorem, $(\mu\vert_{N\cap H})^N$ is a summand of $\mu^G\vert_N=\chi\vert_N=e\theta.$ Hence $(\mu\vert_{N\cap H})^N=f\theta$ for some $f.$ By Theorem K,
this can only happen if $N\cap H=H,$ so $H\supseteq N.$ Then $\mu\vert_N=f\theta$ as claimed.
\end{proof}

In particular, if $N\vartriangleleft G$ is a normal subgroup of odd order and $\theta\in\Irr(N)$ is non-linear, quasi-primitive and
stable in $G,$ then Corollary L shows that no $\chi\in\Irr(G\vert\theta)$ can be monomial. Example 6.4 in \cite{Isaacs1981} shows that
this need not hold when $\theta$ is not stable in $G.$

Our remaining results are direct corollaries of Theorem C, and only use the existence of a subgroup $U$ with
$\bar\chi\chi=(1_U)^G.$ For simplicity, we restrict our attention to solvable groups in the statements of these results.

\begin{corM}
Let $G$ be solvable and suppose $\chi\in\Irr(G)$ is primitive and $\chi(1)$ is odd. Then $\abs{\chi(x)}\leq\abs{\chi(x^m)}$ for any
integer $m.$
\end{corM}

\begin{proof}
By Theorem C, there exists $U\subseteq G$ with $\bar\chi\chi=(1_U)^G.$ Thus for $x\in G,$ $\abs{\chi(x)}^2$ is the number of fixed
points of $x$ in its coset representation on $U.$ But all such points are also fixed by $x^m,$ and the result follows.
\end{proof}

In particular, in Corollary M, if $\chi(x)\neq 0$ then $\chi(x^m)\neq 0$ for any integer $m.$ This improves the
``only if'' side of Theorem A, when $\chi(1)$ is odd.

As further corollaries of Theorem C, we mention some results about zeros of characters of solvable groups. Following
\cite{Gallagher2006}, we write $n_\chi=\mathrm{lcm}\set{o(x)\vert\chi(x)\neq 0},$ the least common multiple of the orders of elements
in the support of $\chi.$ If $G$ is solvable then by \cite[Corollary 4.3]{Wilde2006}, for any $\chi\in\Irr(G),$ $n_\chi$ divides
$\abs{G}/\chi(1).$ When $\chi$ is primitive of odd degree, we have:

\begin{corN}
Let $G$ be solvable and suppose $\chi\in\Irr(G)$ is primitive and $\chi(1)$ is odd. Then $\chi(1)^2n_\chi$ divides $\abs{G}.$
\end{corN}

\begin{proof}
By Theorem C, there exists $U\subseteq G$ with $\bar\chi\chi=(1_U)^G.$ Hence if $x\in G$ and $\chi(x)\neq 0,$ then some conjugate of $x$ lies in $U.$ Hence the order of $x$ divides $\abs{U}=\abs{G}/\chi(1)^2,$ as required.
\end{proof}

\begin{corO}
Let $G$ be solvable and suppose $\chi\in\Irr(G)$ and $\chi(1)$ is odd. Suppose $n_\chi=\abs{G}/\chi(1).$ Then $\chi$ is monomial.
\end{corO}

\begin{proof}
We have $\chi=\psi^G$ for some subgroup $H\subseteq G$ and primitive $\psi\in\Irr(H)$ of odd degree. If $\chi(x)\neq 0$ then for some
$g\in G,$ $x^g$ lies in $H$ and has $\psi(x^g)\neq 0.$ By Corollary N, the order of $x$ divides $\abs{H}/\psi(1)^2.$ Hence $n_\chi$
divides $\abs{H}/\psi(1)^2=\abs{G}/(\chi(1)\psi(1)).$ But $n_\chi=\abs{G}/\chi(1)$ by supposition. Hence $\psi(1)=1$ and $\chi$ is
monomial.
\end{proof}

Two final applications of Theorem C are as follows.

\begin{corP}
Let $G$ be solvable. Suppose $\chi\in\Irr(G)$ has odd degree and let $\pi$ be the set of primes dividing $\chi(1).$ Suppose that
$\chi=\prod\psi_i$ may be written as a product of primitive characters $\psi_i\in\Irr(G).$  Let $H\subseteq G$ and let $H_\pi$ be a
Hall $\pi$-subgroup of $H.$ Then $\chi\vert_H\in\Irr(H)$ if and only if $\chi\vert_{H_\pi}\in\Irr(H_\pi).$
\end{corP}

\begin{proof}
The ``if'' part is clear. For the ``only if'' part, since each $\psi_i(1)$ is odd, then by Theorem C, for each $i$ there exists a
permutation representation $G\times\Omega_i\rightarrow\Omega_i$ whose character is $\bar\psi_i\psi_i.$ Then the action of $G$ on $\Omega=\prod_i\Omega_i$ affords $\bar\chi\chi,$ and since $[\bar\chi\chi,1]=[\chi,\chi]=1,$ it follows that $G$ acts transitively 
on $\Omega$ and we have $\bar\chi\chi=(1_K)^G,$ where $K$ is the stabilizer of a point of $\Omega.$ 
Now if $L\subseteq G$ is any subgroup, then $\chi\vert_L$ is irreducible just when $[\bar\chi\chi,1]_K=[(1_L)^G,1]_K=1,$ that is, just when $LK=G.$ Hence $HK=G,$ and we must show that $H_\pi K=G.$ But $\chi(1)^2=\abs{G:K}$ is a $\pi$-number, so $\abs{G}_{\pi^\prime}$ divides $\abs{K}.$ Hence $\abs{H_\pi
K}=\abs{H}_\pi\abs{K}/\abs{H_\pi\cap K}$ is divisible by $\abs{G}_{\pi^\prime}$ and also by $\abs{H}_\pi\abs{K}_\pi/\abs{H\cap
K}_{\pi}=\abs{HK}_\pi=\abs{G}_\pi.$ Hence $H_\pi K=G$ as required.
\end{proof}

We use Corollary P to generalize \cite[Theorem 6.23]{Isaacs1973}. If $H\subseteq G$ and $\chi\in\Irr(G),$ we say that $\chi$ is a
\emph{relative $M$-character with respect to $H$} if there exists $K\supseteq H$ and $\varphi\in\Irr(K)$ with $\varphi\vert_H\in\Irr(H)$
and $\varphi^G=\chi.$ (When $H$ is a normal subgroup, this is \cite[Definition 6.21]{Isaacs1973}.)

\begin{corQ}
Let $G$ be solvable and suppose $H\subseteq G$ has abelian Sylow $p$-subgroups for each prime $p.$ Let $\chi\in\Irr(G),$ and suppose
that at least one of $\chi(1)$ or $\abs{H}$ is odd. Suppose $\chi$ is a relative $M$-character with respect to $H.$ Then $\chi$ is
monomial.
\end{corQ}

\begin{proof}
By hypothesis, there exist $K\supseteq H$ and $\varphi\in\Irr(K)$ with $\varphi\vert_H\in\Irr(H)$ and $\varphi^G=\chi.$ We may write $\varphi=\mu^K$ for some $L\subseteq K$ and primitive $\mu\in\Irr(L).$ By the theory of special characters (e.g. \cite[Corollary 21.8]{ManzWolf1993}), there exist characters $\mu_p\in\Irr(L)$ with $\mu=\prod_p\mu_p,$ where $p$ runs over the primes dividing $\mu(1).$ Each $\mu_p$ is primitive, since $\mu$ is. Since $\varphi\vert_H=\mu^K\vert_H$ is irreducible, certainly $\mu\vert_{L\cap H}$ is irreducible, and hence so is $\mu_p\vert_{L\cap H}.$ Since at least one of $\mu(1)$ or $\abs{L\cap H}$ is odd, actually $\mu(1)$ must be odd, so Corollary P shows that $\mu_p$ remains irreducible on restriction to the Sylow $p$-subgroups of $L\cap H.$ These, however, are abelian, so $\mu_p$ is linear for each $p.$ Thus $\mu$ is linear and $\chi=\mu^G$ is monomial, as required.  \end{proof}

\begin{rems}
1. Theorem C is false without the condition that $\chi(1)$ be odd. For example, let $G$ be the general linear group $GL_2(\mathbb F_3).$ Then $G$ is a solvable group of order $48,$ and $G$ has an irreducible primitive character $\chi$ of degree $2.$ Also, $G$ has elements $x$ and $y$ of order $3$ and $8$ respectively, such that $\chi(x)$ and $\chi(y)$ are not zero and $\chi(y^2)=0.$ The latter is contrary to the conclusion of Corollary M, and shows that no equation $\bar\chi\chi=(1_U)^G$ can hold. The order of $y$ does not divide $48/\chi(1)^2=12,$ contrary to the conclusion of Corollary N, and $n_\chi=24=\abs{GL_2(\mathbb{F}_3)}/\chi(1),$ but $\chi$ is not monomial contrary to the conclusion of Corollary O.

2. In Corollary O, the condition $n_\chi=\abs{G}/\chi(1)$ is satisfied infinitely often. For example, suppose
$G$ has squarefree order and let $\chi$ be any irreducible character of $G$ of odd degree. If $p$ is a prime dividing
$\abs{G}/\chi(1),$ and $x\in G$ has order $p,$ then $p$ cannot divide $\chi(1),$ and it follows that $\chi(x)\neq 0$ and hence that
$p$ divides $n_\chi.$ Therefore $n_\chi=\abs{G}/\chi(1).$ (Of course, any character of a group of squarefree order is monomial.)
\end{rems}

\section{Proof of Theorems A and B}
To prove Theorems A and B and provide the basis for our other results, we first study \emph{strongly irreducible} characters, as defined by
Brauer in \cite{Brauer1977}. Recall that if $\chi\in\Irr(G),$ then $Z(\chi)$
is the pre-image in $G$ of the centre of $G/\ker\chi.$

\begin{defn}[\cite{Brauer1977}, Definition 1.1]
Let $\chi\in\Irr(G).$ Then $\chi$ is \emph{strongly irreducible} if whenever $N\vartriangleleft G,$ then
either $N\subseteq Z(\chi)$ (so that $\chi\vert_N=\chi(1)\lambda$ for some linear character $\lambda\in\Irr(N)$), or $\chi\vert_N$ is
irreducible.
\end{defn}

Note that a strongly irreducible character is, in particular, quasi-primitive. We also have the following easy observation.

\begin{lem}\label{triviallemma}
Let $\chi\in\Irr(G)$ and let $N\vartriangleleft G$ with $N\subseteq\ker\chi.$ Then $\chi$ is strongly irreducible if and only if $\chi$ is strongly irreducible viewed as a character of $G/N.$
\end{lem}
\begin{proof}
It is immediate from the definition that $\chi\in\Irr(G)$ is strongly irreducible if and only if for any $L\vartriangleleft G,$ either $[\chi,\chi]_L=1$ or $[\chi,\chi]_L=\chi(1)^2.$ However, since $N\subseteq\ker\chi,$ we have $[\chi,\chi]_L=[\chi,\chi]_{LN/N},$ where on the left, we view $\chi\in\Irr(G/N).$ The result follows.
\end{proof}

When $\chi\in\Irr(G)$ is strongly irreducible and $G$
satisfies a suitable solvability hypothesis, then we have restrictive information about $G$ and $\chi,$ as follows.

\begin{prop}\label{prop2}
Let $G$ be $p$-solvable and suppose $\chi\in\Irr(G)$ is faithful and has $\chi(1)$ divisible by $p.$ Let $Z$ be
the centre of $G,$ and let $K\vartriangleleft G$ with $Z\subset K$ be such that $K/Z$ is a minimal normal subgroup of $G/Z.$ Suppose
that $\chi\vert_K$ is irreducible. (In particular, this is the case if $\chi$ is strongly irreducible.) Then the following hold.

\begin{enumerate}
\item $K/Z$ is an elementary abelian $p$-group.
\item For any $x\in G,$ $$\abs{\chi(x)}^2=\begin{cases}\abs{C_K(x):Z} & \text{ if } C_{K/Z}(xZ)=C_K(x)/Z,$$ \\0 & \text{ otherwise.}\end{cases}$$
\item $K/Z=C_{G/Z}(K/Z),$ and is the unique minimal normal subgroup of $G/Z.$
\item $\chi$ is strongly irreducible.
\end{enumerate}
\end{prop}

Parts of Proposition \ref{prop2} are contained in \cite{Brauer1977} (see particularly Proposition 3C and Corollary 3D), but we offer a different approach, based on the following elementary lemma, which is essentially \cite[Problem 3.12]{Isaacs1976}. If $x$ and $y$ are elements of a group $G$ then $[x,y]$ denotes their commutator: $[x,y]=x^{-1}y^{-1}xy.$

\begin{lem}\label{lem1}
Let $G$ be a group and let $\chi\in\Irr(G).$ Suppose $K\vartriangleleft G$ and $\chi\vert_K\in\Irr(K).$ Then for any $x\in G,$
$$\abs{\chi(x)}^2=\frac{\chi(1)}{\abs{K}}\sum_{k\in K}\chi([x,k]).$$
\end{lem}

\begin{proof}
Fix $x\in G$ and let $H=K\langle x\rangle.$ Then $K\subseteq H,$ so certainly $\chi\vert_H\in\Irr(H).$ Let $\mathcal{K}\in Z(\mathbb{C}H)$
be the class sum over $x^H=x^K,$ and let $\mathfrak{R}$ be a representation affording $\chi.$ Then
$\mathfrak{R}(\mathcal{K})=\omega(x)I,$ where $I$ is the identity matrix and $\omega(x)=\vert x^K\vert\chi(x)/\chi(1).$ Hence
$\mathfrak{R}(x^{-1}\mathcal{K})=\omega(x)\mathfrak{R}(x^{-1}).$ Taking traces gives $\omega(x)\chi(x^{-1})=\sum_{y\in
x^K}\chi(x^{-1}y),$ or $$\abs{\chi(x)}^2=\frac{\chi(1)}{\abs{x^K}}\sum_{y\in x^K}\chi(x^{-1}y)=\frac{\chi(1)}{\abs{K}}\sum_{k\in
K}\chi(x^{-1}x^k)$$ as required.
\end{proof}

\begin{proof}[Proof of Proposition \ref{prop2}]
(1) By hypothesis $\chi\vert_K\in\Irr(K),$ so $\chi(1)$ divides $\abs{K/Z(K)}$ and so certainly $p$ divides $\abs{K/Z}.$ Since $K/Z$ is a
chief factor of the $p$-solvable group $G,$ $K/Z$ is an elementary abelian $p$-group as stated.

(2) We have $\chi\vert_Z=\chi(1)\lambda,$ where $\lambda$ is a faithful linear character of $Z.$ Since $\chi(1)>1,$ it follows from
\cite[Theorem 6.18]{Isaacs1976} (we have case (b) of that theorem) that $\chi$ and $\lambda$ are fully ramified with respect to $K/Z.$
This means that $\chi(1)^2=\abs{K:Z}$ and $\chi$ vanishes on $K\backslash Z.$ Hence, Lemma \ref{lem1} gives, for any $x\in G,$ the expression
$$\abs{\chi(x)}^2=\frac{\chi(1)^2}{\abs{K}}\sum_{k\in K\vert [x,k]\in Z}\lambda([x,k])= \frac{1}{\abs{Z}}\sum_{k\in K\vert [x,k]\in
Z}\lambda([x,k]).$$ Write $U_x=C_G(x\mod Z)=\set{k\in K\vert [x,k]\in Z}$ for the range of summation. We have $Z\subseteq
C_K(x)\subseteq U_x$ and $U_x/Z=C_{K/Z}(xZ).$ The map $U_x\rightarrow Z$ given by $k\mapsto [x,k]$ is a homomorphism. Lifting
$\lambda$ to $U_x$ via this map, the formula above becomes $$\abs{\chi(x)}^2=\abs{U_x:Z}[\lambda,1]_{U_x}.$$ The right hand side is
non zero if and only if $\lambda=1_{U_x},$ that is if and only if $U_x=C_K(x),$ and its value in that case is
$\abs{U_x:Z}=\abs{C_K(x):Z}$ as required.

(3) First, we show that $C_G(K)=Z.$ Let $x\in C_G(K).$ Then $K=C_K(x)\subseteq U_x\subseteq K,$ so by part (2),
$\abs{\chi(x)}^2=\abs{K:Z}=\chi(1)^2.$ Since $\chi$ is faithful, it follows that $x\in Z,$ so indeed $C_G(K)=Z.$

Now let $L\supseteq K$ be such that $L/Z=C_{G/Z}(K/Z).$ Since $C_G(K)=Z,$ if $x\in L\backslash Z$ then $C_K(x)\neq K,$ but by
hypothesis $C_{K/Z}(xZ)=K/Z.$ By part (2), $\chi(x)=0.$ Hence $\chi$ vanishes on $L\backslash Z.$ Since $L\supseteq K,$ certainly
$\chi\vert_L\in\Irr(L)$ and $[\chi,\chi]_L=[\chi,\chi]_Z/\abs{L:Z}=1.$ Therefore $\chi(1)^2=[\chi,\chi]_Z=\abs{L:Z}.$ But
$\chi(1)^2=\abs{K:Z},$ so $L=K$ and we conclude that $C_G(K/Z)=K/Z$ as desired. The last statement is clear, since any two distinct
minimal normal subgroups of any group must intersect trivially, and hence must centralize each other.

(4) Let $N\vartriangleleft G$ be any normal subgroup with $N\nsubseteq Z.$ We have to show that $\chi\vert_N$ is irreducible. Now $NZ/Z$ is
a non-trivial normal subgroup of $G/Z,$ so by part (3), $K\subseteq NZ.$ Hence $\chi\vert_{NZ}\in\Irr(NZ)$ and so clearly
$\chi\vert_N\in\Irr(N).$ This completes the proof.
\end{proof}

The next result was stated by Brauer in a more general form (\cite[Proposition 6E]{Brauer1977}) and is also contained in (\cite[Theorem 1.13 and Corollary 1.14]{FergusonTurull1985}). As these authors show, the solvability hypothesis is not actually necessary, but it will not restrict our use of the result.

\begin{prop}\label{prop1}
Suppose $\chi\in\Irr(G)$ is quasi-primitive and assume $G$ is $\pi(\chi(1))$-solvable. Then there exist a central extension $\pi:\tilde G\twoheadrightarrow G$ and strongly irreducible characters $\rho_1,...,\rho_n\in\Irr(\tilde G)$ with $\prod_{i=1}^n\rho_i=\tilde\chi,$ where $\tilde\chi=\chi\circ\pi$ is the lift of $\chi$ to $\tilde G.$
\end{prop}

The proof requires the following lemma, which may also be found at \cite[Lemma 1.1]{FergusonTurull1985} and \cite[Proposition 3C]{Brauer1977}.

\begin{lem}\label{chi1chi2lemma} Let $\chi_1,\chi_2\in\Irr(G)$ and suppose $\chi_1\vert_{\ker\chi_2}$ is irreducible and $\chi_1\chi_2$ is irreducible and quasi-primitive. Then $\chi_2$ is quasi-primitive.
\end{lem}

\begin{proof}
Suppose $\chi_2$ is not quasi-primitive and let $N\vartriangleleft G$ be such that $\chi_2\vert_N$ is not homogeneous. If
$\chi_2\vert_{N\ker\chi_2}$ is homogeneous, then $\chi_2\vert_{N\ker\chi_2}=f\theta$ for some
$\theta\in\Irr(N\ker\chi_2)$ and clearly $\ker\chi_2=\ker\theta,$ so $\theta_N$ is irreducible and $\chi_2\vert_N$ is homogeneous, contradicting the choice of $N.$ Hence we may replace $N$ with $N\ker\chi_2$ and assume that $N\supseteq \ker\chi_2.$
If $\eta$ and $\eta^g$ are distinct constituents of $\chi_2\vert_N$ then in $N,$ $\chi_1\vert_{\ker\chi_2}$ is irreducible and hence
$\chi_1\eta$ and $\chi_1\eta^g=(\chi_1\eta)^g$ are distinct constituents of $\chi_1\chi_2\vert_N.$ This contradiction completes the
proof.
\end{proof}

\begin{proof}[Proof of Proposition \ref{prop1}]
We may assume that $\chi$ is not strongly irreducible. In particular, $\chi$ is non-linear so $Z(\chi)\neq G,$ and we may choose
$K\vartriangleleft G$ be such that $K/Z(\chi)$ is a minimal normal subgroup of $G/Z(\chi).$ If $\chi\vert_K$ is irreducible, then
by Proposition \ref{prop2}(4) applied in $G/\ker\chi,$ and Lemma \ref{triviallemma}, $\chi$ is strongly irreducible, a contradiction.
Hence, since $\chi$ is quasi-primitive, $\chi\vert_K=f\eta$ for some $f>1$ and $\eta\in\Irr(K),$ where $\eta(1)>1$ because $K\nsubseteq
Z(\chi).$ 

By theorems of Clifford (e.g. \cite[Theorem 8.16]{Navarro1998} or \cite[Theorem 51.7]{CurtisReiner1962}) and Schur (e.g. \cite[Theorem 11.17]{Isaacs1976} or \cite[Theorem 53.7]{CurtisReiner1962}), there is a central extension $\pi: \tilde G\twoheadrightarrow G$ such that if $\tilde\chi$ denotes the lift of $\chi$ to $\tilde G,$ then there exist irreducible characters $\chi_1,\chi_2\in\Irr(\tilde G)$ satisfying
$$\chi_1\chi_2=\tilde\chi,$$ $$\chi_1(1)=\eta(1),\chi_1\vert_{\tilde K}\in\Irr(\tilde K),$$ $$\chi_2(1)=e,\tilde K\subseteq Z(\chi_2),$$ where $\tilde K$ is the pre-image of $K$ in $\tilde G.$ Since $\chi\vert_{\ker(\chi_2)}=e\chi_1\vert_{\ker(\chi_2)}$ is homogeneous, we know that $\chi_1\vert_{\ker(\chi_2)}$ is homogeneous, and as $\chi_1\vert_{Z(\chi_2)}$ is irreducible and $Z(\chi_2)/\ker(\chi_2)$ is cyclic, it follows for example by Corollaries 11.22 and 6.17 in \cite{Isaacs1976}, that $\chi_1\vert_{\ker\chi_2}$ is irreducible, so by Lemma \ref{chi1chi2lemma}, $\chi_2$ is quasi-primitive. We show that $\chi_1$ is strongly irreducible. Let $K_1=\tilde KZ(\chi_1).$ Then $K_1\vartriangleleft\tilde G$ and $\chi_1\vert_{K_1}$ is irreducible. Moreover (as for any product of characters), we have $$Z(\tilde\chi)=\set{g\in G_1,\abs{\chi_1(g)}\abs{\chi_2(g)}=\chi_1(1)\chi_2(1)}=Z(\chi_1)\cap Z(\chi_2),$$ and since also $Z(\chi_2)\supseteq \tilde K$ and $Z(\tilde\chi)\subseteq \tilde K,$ it follows that $Z(\tilde\chi)=\tilde K\cap Z(\chi_1),$ and $$K_1/Z(\chi_1)=\tilde K/\tilde K\cap Z(\chi_1)=\tilde K/Z(\tilde\chi)$$ as $\tilde G$-operator groups. Since
$\tilde K/Z(\tilde\chi)$ is a chief factor of $\tilde G,$ we have that $K_1/Z(\chi_1)$ is a minimal normal subgroup of $\tilde
G/Z(\chi_1).$ By Proposition \ref{prop2}(4) applied in $\tilde G/\ker\chi_1,$ and Lemma \ref{triviallemma}, $\chi_1$ is strongly
irreducible. Since $\chi_2(1)<\chi(1),$ an obvious induction applies to give an extension $\pi:G_1\twoheadrightarrow G$ and strongly
irreducible characters $\rho_1,...,\rho_n$ with $$\prod_{i=1}^n\rho_i=\tilde\chi,$$ where $\tilde\chi=\chi\circ\pi$ is the inflation
of $\chi$ to $G_1.$ 
Lastly, we show that we can replace $G_1$ with a central extension of $G.$ Let $A=\bigcap_i\ker\rho_i\cap\ker\pi.$ Then $\pi$ factors through $G_1/A$ and the $\rho_i$ all deflate to $G_1/A.$ Moreover, by Lemma \ref{triviallemma}, each $\rho_i$ is a strongly irreducible character of $G_1/A.$ Hence we may replace $G_1$ with $G_1/A.$ We suppose this replacement has been made, so $A=1.$ Now $$\ker\pi\subseteq\ker\tilde\chi\subseteq Z(\tilde\chi)=\bigcap_{i=1}^n Z(\rho_i),$$ so if $g\in\ker\pi,$ then $[g,\tilde G]\subseteq \bigcap_{i=1}^n\ker\rho\cap\ker\pi=A=1,$ and so $\ker\pi$ is central. This completes the proof.
\end{proof}

\begin{proof}[Proof of Theorem B]
By Proposition \ref{prop1} there exist a central extension $\pi:\tilde G\rightarrow G$ and strongly irreducible characters $\rho_i\in\Irr(G)$ with $\prod_i\rho_i=\tilde\chi,$ where $\tilde\chi=\chi\circ\pi.$ Let $g\in G$ and let $g_1\in\tilde G$ with $\pi(g_1)=g.$ Then $\abs{\chi(g)}^2=\prod_i\abs{\rho_i(g_1)}^2.$ Hence, we need only prove the result each $\rho_i.$ If $\rho_i$ is linear then this is clear. Otherwise, let $p$ be a prime dividing $\rho_i(1).$ Then by hypothesis $G$ is $p$-solvable and hence so is $\tilde G,$ being a central extension of $G.$ The result follows from Lemma \ref{triviallemma} and Proposition \ref{prop2}(2) applied in $\tilde G/\ker\rho_i.$
\end{proof}

For the proof of Theorem A, we need the following additional information.

\begin{prop}\label{prop3}
Let $G$ be $p$-solvable and suppose that $\chi\in\Irr(G)$ is strongly irreducible and $p$ divides $\chi(1).$ Let $x\in G.$ Then
$\chi(x)\neq 0$ if and only if $\chi(x_p)\neq 0.$
\end{prop}

\begin{proof}
Let $N=\ker\chi.$ By Lemma \ref{triviallemma}, $\chi$ is strongly irreducible viewed as a character of $G/N.$ Since $(xN)_p=x_pN,$ the result will follow if we can show that $\chi(xN)\neq 0$ if and only if $\chi((xN)_p)\neq 0.$ Thus we may assume $\chi$ is faithful, and we are in the
situation of Proposition \ref{prop2}. We use the notation of that Proposition. Thus, we must show that $U_x=C_K(x)$ if and only if $U_{x_p}=C_K(x_p).$

If $k\in U_x,$ then $k^x=kz$ for some $z\in Z.$ Since $K/Z$ is elementary abelian, $k^p\in Z.$ Hence $(k^x)^p=k^p=k^pz^p$ and $z^p=1.$
Then since $kxk^{-1}=xz$ and $z$ is a $p$-element, we have $kx_{p^\prime}k^{-1}=x_{p^\prime}$ and $kx_pk^{-1}=x_pz.$ Hence $U_x=U_{x_p}\cap C_K(x_{p^\prime}).$ Thus if $U_{x_p}=C_K(x_p),$ then $U_x=C_K(x).$

On the other hand, suppose $U_{x_p}\neq C_K(x_p),$ and choose $k$ in $U_{x_p}\backslash C_K(x_p).$ Thus $k^{x_p}=kz$ for some $z\in Z$ with $z\neq 1.$ As above, since $k^p\in Z,$ we have $z^p=1.$ Let $$k^\prime=\prod_{r=1}^{o(x_{p^\prime})}k^{((x_{p^\prime})^r)}=k^{x_{p^\prime}}k^{(x_{p^\prime})^2}...k.$$ Since $K/Z$ is abelian, $(k^\prime)^{x_{p^\prime}}\in k^\prime Z.$ Also, 
$$(k^\prime)^{x_p}=\prod_{r=1}^{o(x_{p^\prime})}(kz)^{((x_{p^\prime})^r)}=k^\prime z^{o(x_{p^\prime})}.$$ Therefore $k^\prime\in U_{x_p}\cap U_{x_{p^\prime}}=U_x,$ but since $z\neq 1$ and $z^p=1,$ we have $z^{o(x_{p^\prime})}\neq 1,$ so $k^\prime\notin C_K(x_p).$ Thus $k^\prime\in
U_x\backslash C_K(x),$ and it follows that $U_x\neq C_K(x).$ This completes the proof.
\end{proof}

\begin{proof}[Proof of Theorem A]
As in the proof of Theorem B above, we have by Proposition \ref{prop1} a central extension $\pi:\tilde G\rightarrow G$ and strongly irreducible characters $\rho_i\in\Irr(G)$ with $\prod_i\rho_i=\tilde\chi,$ where $\tilde\chi=\chi\circ\pi.$ For $g\in G,$ let $g_1\in\tilde G$ with $\pi(g_1)=g.$ Then $\pi(g_{1p})=g_p.$ We have $\chi(g)\neq 0$ if and only if $\rho_i(g_1)\neq 0$ for each $i,$ which by Proposition \ref{prop3} is equivalent to $\tilde\chi(g_{1p})=\chi(g_p)\neq 0,$ as required.
\end{proof}

\section{Properties of solvable groups}
Before addressing Theorems C and D in the next section, we collect here some mostly well known group-theoretic and cohomological facts about
solvable groups and their maximal subgroups. We also give the proofs of Theorems G and H.

\begin{prop}\label{CVprop}
Let $G$ be $p$-solvable and let $X$ be a faithful, simple $G$-module with $p\vert\abs{X}.$ Then $H^1(G,X)=H^2(G,X)=0.$
\end{prop}

\begin{proof}
Since $X$ is simple, it is an elementary abelian $p$-group and hence is a faithful, simple $\mathbb{F}_pG$-module. As
$X_{\mathbf{O}_p(G)}$ is semi-simple, $\mathbf{O}_p(G)$ acts trivially on $X$ and as $X$ is faithful, $\mathbf{O}_p(G)=1.$ Then since
$G$ is $p$-solvable, $\mathbf O_{p^\prime}(G)\neq 1.$ Let $L=\mathbf O_{p^\prime}(G).$ As $X$ is faithful and simple,
$H^0(L,X)=X^{L}=0.$ Also, $H^q(L,X)=0$ for all $q>0$ because $\abs{L}$ is prime to $p.$ In the Lyndon-Hochschild-Serre spectral sequence
$E_2^{p,q}=H^p(G/L,H^q(L,X))\Rightarrow H^{p+q}(G,X),$ $E_2^{p,q}=0$ for all $p,q.$ Thus $H^*(G,X)=0,$ and in particular
$H^1(G,X)=H^2(G,X)=0$ as desired.
\end{proof}

We record two well known corollaries of this result.

\begin{cor}\label{maxcor1}
Let $G$ be $p$-solvable and suppose $X$ is a minimal normal subgroup of $G$ with $p$ dividing $\abs{X}.$ Suppose
$X=C_G(X).$ Then $X$ has a unique conjugacy class of complements in $G.$ These complements are maximal subgroups of $G.$
\end{cor}

\begin{proof}
Clearly $X$ is elementary abelian of order a power of $p.$ View $X$ as a $G/X$-module. Then $X$ is faithful and simple, so by
Proposition \ref{CVprop}, $H^1(G/X,X)=H^2(G/X,X)=0.$ By the usual interpretations of the first and second cohomology groups (see
e.g. \cite{Brown1982}, Theorem 3.12 and Proposition 2.1), $X$ has a unique conjugacy class of complements in $G.$ Let $U$ be one of
these complements. It remains to show that $U$ is maximal in $G.$ Suppose $U\subseteq L\subseteq G$ with $U\neq L.$ Then $L=U(L\cap
N)$ by the modular law, so $L\cap X\neq 1,$ and since $X$ is abelian, $L\cap X\vartriangleleft XL=G.$ Since $X$ is minimal normal,
$L\subseteq X$ and so $L=G.$ Thus $U$ is maximal as required.
\end{proof}

\begin{cor}\label{firstcohomologylemma}
Let $X$ be a simple module for the group $G,$ and suppose $G$ is $p$-solvable where $p$ is the characteristic of $X.$ Then restriction
is an isomorphism $H^1(G,X)\xrightarrow{\res}\hm_G(C_G(X),X).$
\end{cor}
\begin{proof}
We use the five-term exact sequence (\cite[Corollary 7.2.3]{Evens1991}),
$$1\rightarrow H^1(G/C_G(X),X)\xrightarrow{\inf} H^1(G,X)\xrightarrow{\res} H^1(C_G(X),X)^G\xrightarrow{d_2}H^2(G/C_G(X),X).$$
Since $X$ is a faithful, simple module for $G/C_G(X),$ we have $H^{1,2}(G/C_G(X),X)=0$ by Proposition \ref{CVprop}. Therefore the middle restriction map in this sequence is an isomorphism $H^1(G,X)\rightarrow H^1(C_G(X),X)^G.$ Since $C_G(X)$ acts
trivially on $X,$ we have $H^1(C_G(X),X)^G=\hm_G(C_G(X),X),$ and the proof is complete.
\end{proof}

Recall that the Schur multiplier of a group $G$ is the cohomology group $M(G)=H^2(G,\mathbb{C}^\times).$ Proposition \ref{prop2} shows
that when $\chi\in\Irr(G)$ is non-linear and strongly irreducible, where $G$ is $p$-solvable for the prime dividing $\chi(1),$ then
$G/Z(G)$ has a unique, minimal normal subgroup which is abelian and self-centralizing. In this situation, the structure of the Schur
multiplier of $G$ is known, and is described in the next proposition.

\begin{prop}\label{MGprop}
Suppose $G$ is $p$-solvable and has a self-centralizing minimal normal subgroup $X=C_G(X),$ where $X$ is an abelian $p$-group. Let $U$
be a complement to $X$ in $G$ (such a complement exists by Corollary \ref{maxcor1}). Then
$$M(G)\cong M(U)\times M(X)^G,$$ where the isomorphism is given by the pair of maps $(\res_U,\res_X).$
\end{prop}

\begin{proof} 
Consider the sequence of groups and natural maps $$1\rightarrow M(G/X)\xrightarrow{\mathrm{inf}_{G/X}}M(G)\xrightarrow{\mathrm{res_X}}M(X)^G\rightarrow 1.$$ 
The complement $U$ gives rise to a splitting of this sequence. Thus if $\iota:G/X\rightarrow U$ is the natural isomorphism which takes $uX\mapsto u$ for $u\in U,$ then $\iota^*$ is an isomorphism $M(U)\rightarrow M(G/X)$ and the map $\iota^*\res_{U}:M(G)\rightarrow M(G/X)$ satisfies $\iota^*\res_{U}\mathrm{inf}_{G/X}=1_{M(G/X)}.$ This shows that the sequence is split and is also exact at $M(G/X).$ Let $\tilde M(G)\subseteq M(G)$ be the kernel of restriction to $U.$ Then $M(G)\cong M(U)\times \tilde M(G)$ as an internal direct product, where we identify $M(U)\cong M(G/X)$ with its image under inflation. By a theorem of Tahara (\cite[Theorem 2(II)]{Tahara1972}, or
see \cite[Theorem 2.2.5]{Karpilovsky1987}), we have an exact sequence 
$$1\rightarrow H^1(U,X^*)\rightarrow \tilde M(G)\xrightarrow{\res_X}M(X)^G\rightarrow H^2(U,X^*).$$ By Proposition \ref{CVprop},
since $U$ is $p$-solvable and $X^*$ is a faithful, simple module for $U,$ we have $H^1(U,X^*)=H^2(U,X^*)=0.$ Therefore in this case, restriction to $X$ gives an isomorphism $\tilde M(G)\xrightarrow{\res_X}M(X)^G.$ The result now follows.
\end{proof}

\begin{rem}
Proposition \ref{MGprop} can also be obtained by a spectral sequence argument. We sketch this alternative proof. Relevant
properties of the Lyndon-Hochschild-Serre spectral sequence can be found for example in \cite[Section 7.2]{Evens1991}. Consider again the sequence $$1\rightarrow M(G/X)\xrightarrow{\mathrm{inf}_{G/X}}M(G)\xrightarrow{\mathrm{res_X}}M(X)^G\rightarrow 1.$$ As argued above, the sequence is split and exact at $M(G/X).$ It is easy to see that Proposition \ref{MGprop} is equivalent to the sequence being exact at $M(G)$ and $M(X)^G$ as well, which we prove by considering the Lyndon-Hochschild-Serre spectral sequence $E_{p,q}^2=H^p(G/X,H^q(X,\mathbb C^\times))\Rightarrow H^{p+q}(G,\mathbb C^\times).$
First, since the inflation map $M(G/X)\rightarrow M(G)$ is a monomorphism, it follows that all differentials into the bottom rows $E_{p,0}^r$ of the spectral sequence are zero. By Proposition \ref{CVprop}, $E_2^{p,1}=H^p(G/X,\hm(X,\mathbb C^\times))=0$ for $p\geq 0.$ The differential $d_2^{0,2}:E_2^{0,2}\rightarrow E_2^{2,1}$ is zero, since its image $E_2^{2,1}$ is zero. Also $d_3^{0,2}$ is zero since its image $E_3^{3,0}$ is in the bottom row of its page of the sequence. There are no differentials out of $E_r^{0,2}$ for $r\geq 3,$ so in fact all differentials out of $E_*^{0,2}$ vanish.
The vertical edge homomorphism $H^2(G,\mathbb C^\times)\rightarrow E_\infty^{0,q}\rightarrow E_2^{0,q}=H^2(X,\mathbb C^\times)^G,$
where the first map is projection and the second map is inclusion, is just the restriction map, and we have just shown that the second
map is onto, as there are no non-zero differentials out of $E_*^{0,2}.$ This gives exactness at $M(X)^G.$ Also, since $E_\infty^{1,1}=E_2^{1,1}=0,$ it follows that $\ker\res_X=E_\infty^{2,0}=\mathrm{im}\inf_{G/X},$ which is exactness at $M(G),$ as required.
\end{rem}

\begin{cor}\label{sigmacor}
In the situation of Proposition \ref{MGprop}, fix a particular complement $U$ to $X,$ and define $\sigma:G\rightarrow G$ by
$\sigma(xu)=x^{-1}u,$ where $x\in X$ and $u\in U.$ Then $\sigma\in\Aut(G)$ and $\sigma$ acts trivially on $M(G).$
\end{cor}

\begin{proof}
We have
$$\sigma(x_1u_1)\sigma(x_2u_2)=(x_1^{-1}(u_1x_2^{-1}u_1^{-1})u_1u_2)=((x_1u_1x_2u_1^{-1})^{-1}u_1u_2)=\sigma(x_1u_1x_2u_2),$$
where in the second equality, we have used the fact that $X$ is abelian. Thus $\sigma\in\Aut(G).$ By Proposition \ref{MGprop},
$M(G)\cong M(U)\times M(X)^G,$ and clearly $\sigma$ preserves this direct sum decomposition and acts trivially on $M(U).$ We are
reduced to showing that $\sigma$ acts trivially on $M(X)^G.$ However, since $X$ is elementary abelian, it follows from \cite[Problem
11.16]{Isaacs1976} that each class of $M(X)$ contains a bilinear form $\eta: X\times X\rightarrow\mathbb{C}^\times,$ and we have
$\eta^\sigma(x,y)=\eta(x^\sigma,y^\sigma)=\eta(x^{-1},y^{-1})=\eta(x,y).$
\end{proof}

Next, we recall some properties of maximal subgroups of $p$-solvable groups. We will need the concepts of covering and avoidance. Recall that if $K/L$ is a chief factor of $G$ (that is, $K,L\vartriangleleft G,$ $L\subseteq K,$ and $K/L$ is a minimal normal subgroup of $G/L$), and $U\subseteq G$ is a subgroup, then $U$ \emph{covers} $K/L$ if $K=L(U\cap K),$ and $U$ \emph{avoids} $K/L$ if $U\cap K\subseteq L.$

\begin{lem}\label{calem}
$U$ covers $K/L$ if and only if $UK=UL,$ and $U$ avoids $K/L$ if and only if $\abs{UK:UL}=\abs{K:L}.$
\end{lem}
\begin{proof}
Clearly $UL\subseteq UK.$ If $U$ covers $K/L$ then $UK=U(U\cap K)L=UL,$ and if $UK=UL$ then $K=K\cap UL=L(U\cap K)$ and $U$
covers $K/L.$ For the statements about avoidance, we use $\abs{UK:UL}=\abs{K:L}\abs{U\cap L}/\abs{U\cap K}.$ Clearly $U$ avoids $K/L$ if
and only if $K\cap U=L\cap U,$ and the result follows.
\end{proof}

Let $M$ be a maximal subgroup of $G$ and suppose $G$ is $p$-solvable for a prime dividing $\abs{G:M}.$ In Section 2, we said that we could associate to $M$ a right $\mathbb F_pG$-module $X_M,$ where $X_M$ is the unique minimal normal subgroup of $G/\core_G(M).$ This is justified in the following lemma.

\begin{lem}\label{XMlem}
Let $M$ be a maximal subgroup of $G,$ and suppose $G$ is $p$-solvable for a prime $p$ dividing $\abs{G:M}.$ Then $G/\core_G(M)$ has a unique minimal normal subgroup $X_M,$ and $X_M$ is an elementary abelian $p$-group and satisfies $X_M=C_G(X_M)/\core_G(M).$ Also, $MC_G(X_M)=G$ and $M\cap C_G(X_M)=\core_G(M).$ In particular, $\abs{G:M}=\abs{X_M}$ is a power of $p.$ 
\end{lem}

\begin{proof}
Write $N=\core_G(M)$ and let $X$ be any minimal normal subgroup of $G/N.$ Then $X\nsubseteq L/N,$ so $G/N=X(M/N).$ We have $\abs{X:X\cap M/N}=\abs{G:M},$ so by hypothesis, $G$ is $p$-solvable for a prime dividing $\abs{X},$ and it follows that $X$ is an elementary abelian $p$-group. Therefore $C_{G/N}(X)=XC_{M/N}(X)$ by the modular law. Since $X$ and $M/N$ both normalize $C_{M/N}(X),$ we have $C_{M/N}(X)\vartriangleleft G/N,$ so as $M/N$ is core free, we find that $C_{G/N}(X)=X.$ If $Y$ is a minimal normal subgroup of $G/N$ with $Y\neq X,$ then $[X,Y]\subseteq X\cap Y=1,$ so $Y\subseteq C_{G/N}(X)=X,$ a contradiction. Hence $X$ is unique. We may therefore unambiguously label it $X_M.$ We have already shown that $X_M=C_{G/N}(X_M)=C_G(X_M)/N,$ and the remaining statements are clear.
\end{proof}

\begin{lem}\label{avoidlemma}
Let $M$ be a maximal subgroup of $G,$ and assume $G$ is $p$-solvable for a prime $p$ dividing $\abs{G:M}.$ Let $K/L$ be a chief factor of $G.$ Then exactly one of the following holds.
\begin{enumerate}
\item $M$ covers $K/L.$ 
\item $M$ avoids $K/L,$ $K/L$ is abelian and $K/L\cong X_M$ as $G$-modules. 
\end{enumerate}
Furthermore, in any given chief series of $G,$ $M$ avoids exactly one factor, and covers the other factors.
\end{lem}

\begin{proof}
Assume $M$ does not cover $K/L.$ Then by Lemma \ref{calem}, $KM\neq LM.$ Since $M$ is maximal, $L\subseteq M$ and $KM=G.$ Therefore $\abs{K/L}$ is divisible by $\abs{K:K\cap M}=\abs{G:M}.$ Since $G$ is $p$-solvable for a prime $p$ dividing $\abs{G:M},$ it follows that $K/L$ is an abelian $p$-chief factor of $G.$ Therefore since $L\subseteq K\cap M\subset K,$ we have $K\cap M\vartriangleleft KM=G,$ so $M\cap K=L$ and $M$ avoids $K/L$ as claimed. Also, $L\subseteq K\cap\core_G(M)\subseteq K\cap M=L,$ so $L=K\cap\core_G(M).$ Therefore, as $G$-operator groups, $K/L\cong K\core_G(M)/\core_G(M).$ Hence $K\core_G(M)/\core_G(M)$ is a minimal normal subgroup of $G/\core_G(M),$ so $K\core_G(M)/\core_G(M)=X_M$ by Lemma \ref{XMlem}, and $K/L\cong X_M.$ Hence all the statements of (2) hold. Finally, in any chief series, exactly one factor $K/L$ satisfies $L\subseteq M$ and $KM=G.$ By the first part of the proof, this one chief factor is avoided, and the other factors in the series are covered.
\end{proof}

If $L$ is a subgroup of the group $G,$ we write $L^G=\set{L^g\vert g\in G},$ the set of subgroups of $G$ which are conjugate to $L.$

\begin{lem}\label{lem3}
Let $L,M$ be maximal subgroups of $G,$ and suppose $G$ is $p$-solvable for some prime $p$ which divides either $\abs{G:L}$ or $\abs{G:M}.$ Then $\core_G(L)=\core_G(M)$ if and only if $L^G=M^G.$
\end{lem}

\begin{proof}
The ``if'' part is trivial. On the other hand, if $\core_G(L)=\core_G(M)=N,$ then as in the proof of Lemma \ref{XMlem}, $G/N$ has a unique minimal normal subgroup $X,$ which satisfies $X=C_{G/N}(X).$ $L/N$ and $M/N$ both complement $X$ in $G/N,$ so they are conjugate by Corollary \ref{maxcor1}.
\end{proof}

\begin{prop}\label{goingdownprop}
Let $H$ and $M$ be subgroups of $G,$ and suppose that $M$ is maximal and $HM=G.$ Assume that $G$ is $p$-solvable for a prime $p$ dividing $\abs{G:M},$ (so that $X_M$ is defined). Assume further that $X_M\vert_H$ is
a semi-simple $H$-module, and write $X_M\vert_H=\prod_{i=1}^nX_i,$ where each $X_i$ is a simple $H$-module. Then there exist maximal
subgroups $M_i$ of $H,$ for $1\leq i\leq n,$ such that $$M\cap H=\bigcap_{i=1}^n M_i \text{ and for each $i,$ }X_{M_i}\cong X_i \text{
as $H$-modules}.$$
\end{prop}

\begin{proof}
We use induction on the index $\abs{G:H}.$ We will say that the proposition holds for $(G_1,H_1,M_1)$ if it holds with $G,H,M$ in the statement replaced by $G_1,H_1$ and $M_1$ respectively. Let $H_1=HC_G(X_M),$ so that $H\subseteq H_1\subseteq G,$ and suppose first that both of
these inclusions are proper. $H_1$ induces the same transformations as $H$ on $X_M,$ so we may regard each $X_i$ as an $H_1$-module, and we then have $X_M\vert_{H_1}=\prod_{i=1}^nX_i,$ where $X_i\vert_H$ is simple for each $i.$ In particular, $X_{H_1}$ is semi-simple. Also $H_1M=G.$ Since $\abs{G:H_1}<\abs{G:H},$ the proposition holds for $(G,H_1,M)$ by induction. Thus there exist maximal subgroups $L_1,...,L_n$ of $H_1$ such
that $$M\cap H_1=\bigcap_{i=1}^nL_i \text{ and for each $i,$ }X_{L_i}\cong X_i\text{ as $H_1$-modules.}$$ We have $H_1=H_1\cap
MH=(M\cap H_1)H,$ so for each $i,$ since $L_i\supseteq M\cap H_1,$ we have $H_1=L_iH.$ Also $X_{L_i}\vert_H\cong X_i\vert_H$ is simple. Since
$\abs{H_1:H}<\abs{G:H}$ and $H_1$ is $\pi(\abs{H_1:L_i})$-solvable, by induction the proposition holds for $(H_1,H,L_i)$ for each $i.$ Since $X_i\vert_H$ is simple, the proposition here says that $M_i=L_i\cap H$ is a maximal subgroup of $H$ with $X_{M_i}\cong X_i$ as $H$-modules. Also $M\cap H=(M\cap H_1)\cap H=\bigcap_i M_i,$ and this completes the proof when $H\neq HC_G(X_M)\neq G.$ Hence, we may now assume that one of the following is true.
\begin{enumerate}
\item $C_G(X_M)\subseteq H.$
\item $HC_G(X_M)=G.$
\end{enumerate}
We treat these two cases separately. First, suppose that (1) holds. By Lemma \ref{XMlem},
$$(M\cap H)\cap C_G(X_M)=H\cap (M\cap C_G(X_M))=\core_G(M).$$ Suppose $A,B\vartriangleleft H$ with $\core_G(M)\subseteq A,B\subseteq
C_G(X_M).$ Then $A_1=A(M\cap H)$ and $B_1=B(M\cap H)$ are subgroups of $H.$ If $g=am_1=bm_2\in A_1\cap B_1,$ where $a\in A,$ $b\in B$ and $m_{1,2}\in H\cap M,$ then $$a=bm_2m_1^{-1}\in A\cap B_1\subseteq B_1\cap C_G(X_M)=B(M\cap H\cap C_G(X_M))=B\core_G(M)=B.$$ Hence $A_1\cap B_1=(A\cap B)(M\cap H).$
Now since $X_M=C_G(X_M)/\core_G(M),$ there exist $R_i\vartriangleleft H$ with $\core_G(M)\subset R_i\subset C_G(X_M),$ $\core_G(M)/R_i\cong X_i,$ and $$C_G(X_M)/\core_G(M)=C_G(X_M)/R_1\times...\times C_G(X_M)/R_n.$$ Set $M_i=(M\cap H)R_i.$ By the reasoning above, $$\bigcap_{i=1}^nM_i=(\bigcap_{i=1}^nR_i)(M\cap H)=M\cap H.$$ Also, $$M_i\cap C_G(X_M)=(M\cap H)R_i\cap C_G(X_M)=(M\cap H\cap C_G(X_M))R_i=R_i$$ and $$H\supseteq M_iC_G(X_M)\supseteq(M\cap H)C_G(X_M)=H\cap MC_G(X_M)=H,$$ since $MC_G(X_M)=G.$ Thus $M_i$ complements the abelian chief factor $C_G(X_M)/R_i$ of $H.$ So $M_i$ is a maximal subgroup of $H,$ and by Lemma \ref{avoidlemma}, $X_{M_i}\cong X_i$ for each $i.$ This treats case (1).

Next, suppose that (2) holds. Let $L$ be a maximal subgroup of $G$ containing $H.$ Then $LM=G$ and $LC_G(X_M)=G,$ so $X_M\vert_L$ is
simple. If $L\neq H,$ the proposition holds for $(G,L,M)$ by induction, and says that $L\cap M$ is maximal in $L$ with $X_{L\cap M}\cong X_M\vert_L$ as $L$-modules. We have $\abs{L:H}<\abs{G:H}$ and $L$ is $\pi(\abs{L:L\cap M})$-solvable. Hence the proposition also holds for $(L,H,L\cap M).$ Thus $H\cap (L\cap M)=H\cap M$ is maximal in $H$ and has $X_{H\cap M}\cong X_{L\cap M}\vert_H\cong X_M\vert_H,$ which gives the proposition for $(G,H,M)$ in this case.

Hence, we may finally assume that (2) holds and in addition, $H$ is a maximal subgroup of $G.$ We must show that $H\cap M$ is maximal in $H$ and has $X_{H\cap M}\cong X_M\vert_H.$ If $\core_G(M)\subseteq \core_G(H),$ then since $X_M=C_G(X_M)/\core_G(M)$ is the unique minimal normal subgroup of $G/\core_G(M)$ and $C_G(X_M)\nsubseteq H,$ we have $\core_G(H)=\core_G(M)$ and by Lemma \ref{lem3}, $H$ and $M$ are conjugate. Since this is not consistent with $HM=G,$ we conclude that $\core_G(M)\nsubseteq H$ and therefore since $H$ is maximal, $H\core_G(M)=G.$ Then $(H\cap C_G(X_M))\core_G(M)=H\core_G(M)\cap C_G(X_M)=C_G(X_M),$ so as $H$-operator groups, $$\frac{H\cap C_G(X_M)}{H\cap\core_G(M)}\cong\frac{(H\cap C_G(X_M))\core_G(M)}{\core_G(M)}=\frac{C_G(X_M)}{\core_G(M)}.$$ Thus $H\cap C_G(X_M)/H\cap\core_G(M)$ is an abelian chief factor of $H$ which is $H$-isomorphic to $X_M\vert_H.$ Furthermore, since $H\core_G(M)=G,$ we have $M=(M\cap H)\core_G(M)$ and so $G=(M\cap H)C_G(X_M).$ Thus $H=(M\cap H)(H\cap C_G(X_M)).$ Also $(M\cap H)\cap (H\cap C_G(M))=H\cap\core_G(M).$ Therefore $M\cap H$ complements the abelian chief factor $H\cap C_G(X_M)/H\cap\core_G(M)$ of $H,$ so $M\cap H$ is maximal in $H,$ and $X_{M\cap H}\cong X_M\vert_H$ by Lemma \ref{avoidlemma}. This completes the proof.
\end{proof}

If $L\subseteq G$ and $M\subseteq G$ are maximal and $G$ is $p$-solvable for each prime dividing $\abs{G:L}$ or $\abs{G:M},$ then we write $L^G\leq M^G$ if $\core_G(L)\subseteq \core_G(M).$ By Lemma \ref{lem3}, if $L^G\geq M^G$ and $M^G\leq L^G$ then $L^G=M^G,$ so $\geq$ is a partial order on the set of conjugacy classes of such maximal subgroups. The following is well known.

\begin{cor}\label{maxsubgroupcor}
Let $L,M$ be maximal subgroups of $G,$ and assume that $G$ is $p$-solvable for each prime $p$ dividing $\abs{G:L}$ or $\abs{G:M}.$ If $L^G\ngeq M^G,$ then $L\cap M$ is a maximal subgroup of $L.$
\end{cor}
\begin{proof}
By hypothesis, $L\nsupseteq \core_G(M).$ Therefore $L\core_G(M)=G,$ so certainly $LC_G(X_M)=G$ and it follows that $X_M\vert_L$ is a
simple $L$-module. Also $LM=G,$ and the result follows from Proposition \ref{goingdownprop}.
\end{proof}

\begin{lem}\label{greenlemma}
Let $G$ be any group and let $N\vartriangleleft G$ be a normal subgroup of index a power of the prime $p.$ Let $X_i,1\leq i\leq n$ be
pairwise non-isomorphic simple $\mathbb F_pG$-modules. Then $X_i\vert_N=\prod_{j=1}^{m_i}X_{i,j}$ for simple modules $X_{i,j},$ and
the $X_{i,j}$ and $X_{i^\prime,j^\prime}$ are non-isomorphic unless $i=i^\prime$ and $j=j^\prime.$
\end{lem}

\begin{proof} Since $N\vartriangleleft G,$ we certainly have $X_i\vert_N=\prod_{j=1}^{m_i}X_{i,j}$ for simple $\mathbb F_pN$-modules $X_{i,j}.$ We must establish that these modules are pairwise non-isomorphic. Suppose $X_{i,j}\cong X_{i^\prime,j^\prime}$ for some $i,j$ and $i^\prime,j^\prime.$ Then $X_i\vert_N$ and $X_{i^\prime}\vert_N$ have an isomorphic direct summand, so $\hm_{\mathbb F_pN}(X_i,X_j)\neq 0.$ Since $\hm_{\mathbb F_pN}(X_i,X_j)$ has order a power of $p,$ and $G/N$ is a $p$-group, the conjugation action of $G/N$ on $\hm_{\mathbb F_pN}(X_i,X_j)$ must fix some non zero element $\iota\in\hm_{\mathbb F_pN}(X_i,X_j).$ Then $\iota$ is an $\mathbb F_pG$-isomorphism between $X_i$ and $X_{i^\prime},$ so $i=i^\prime.$ To complete the proof, we need to show that for any simple $\hm_{\mathbb F_pG}$-module $X,$ then $X\vert_N$ contains each simple summand with multiplicity at most one. By Clifford theory (written multiplicatively), we have $$X\vert_N=V_1\times ...\times V_n$$ where the $V_i$ are the homogeneous components of $X\vert_N,$ and if $T_i=\set{g\in G\vert V_i^g=V_i}$ then $V_i$ is a simple $\mathbb{F}_pT$ module. But then $V_i$ is simple, by e.g. \cite[Proposition 8.3]{DoerkHawkes1992}, and the proof is complete.
\end{proof}

\begin{proof}[Proof of Theorem H]
We have $G$ $p$-solvable $N\vartriangleleft G$ with index a power of $p,$ $U\subseteq G$ is a regular intersection such that
$\abs{G:U}$ is also a power of $p,$ and $NU=G.$ By hypothesis, there exist maximal subgroups $M_i$ with $U=\bigcap_iM_i$ and the
$X_{M_i}$ pairwise non-isomorphic. Since $N$ is normal, $X_{M_i}\vert_N$ is semisimple and since $NM_i=G,$ Proposition \ref{goingdownprop}
applies, and each $M_i\cap N$ is an intersection of maximal subgroups $L_{i,j}$ of $N$ with $X_i\vert_N\cong\prod_jX_{L_{i,j}}.$ Since
the $X_{M_i}$ are pairwise non-isomorphic, it follows by Lemma \ref{greenlemma} that the $X_{L_{i,j}}$ are non-isomorphic for distinct
$i$ or $j,$ so by definition, $U\cap N=\bigcap_i(M_i\cap N)$ is a regular intersection.
\end{proof}

We end this section with two omnibus lemmas describing the behaviour of complete and regular intersections, and the proof of Theorem G.

\begin{lem}\label{cilemma}
Let $G$ be a group, and suppose $U=\bigcap_{i=1}^nM_i$ is the complete intersection of the maximal subgroups $M_i\subset G,1\leq i\leq n.$ Then
\begin{enumerate}
\item For any $g_1,...,g_n\in G,$ if $U_1=\bigcap_iM_i^{g_i}$ then $U_1$ is $G$-conjugate to $U.$
\item For any $I\subseteq\set{1,...,n},$ $U_I=\bigcap_{i\in I}M_i$ is a complete intersection, and if $J\subseteq\set{1,...,n}$ and $I$ and $J$ are disjoint, then $U_IU_J=G.$
\item If $V\subseteq G$ is also a complete intersection and $UV=G,$ then $U\cap V$ is a complete intersection.
\item If $\pi:G_1\twoheadrightarrow G$ is a central extension, then $\pi^{-1}(U)$ is a complete intersection in $G_1.$ Also, if $W\subseteq G_1$ is a complete intersection and $\ker\pi\subseteq W,$ then $\pi(W)$ is a complete intersection in $G.$
\end{enumerate}
\end{lem}

\begin{proof}
Let $\Omega_i$ be the set of right cosets of $M_i$ with its natural right $G$-action, and let $\Omega=\prod_i\Omega_i.$ Then for any $g_1,...,g_n\in G,$ $(M_1g_1,...,M_ng_n)\in\Omega$ and $$\mathrm{Stab}_G(M_1g_1,...,M_ng_n)=\bigcap_iM_i^{g_i}.$$ Thus the $M_i$ form a complete intersection exactly when $G$ is $\pi(\abs{\Omega})$-solvable, and $\abs{G:\mathrm{Stab}_G(M_1,...,M_n)}=\abs{\Omega},$ or equivalently, when $G$ is $\pi(\abs{\Omega})$-solvable and is transitive on $\Omega.$ If so, then the stabilizers of all points of $\Omega$ are conjugate, which proves (1). For (2), write $\Omega_I=\prod_{i\in I}\Omega_i.$ If $G$ is transitive on $\Omega,$ then $G$ is clearly transitive on $\Omega_I,$ which as above, is equivalent to $U_I=\bigcap_{i\in I}M_i$ being a complete intersection. If $I$ and $J$ are disjoint then again $G$ is transitive on $\Omega_{I\cup J},$ so $U_I\cap U_J$ is the complete intersection of the $M_i$ over $I\cup J,$ and it follows by computing orders that $U_IU_J=G$ as required. (3) similarly follows from the equation $\abs{G:U\cap V}=\abs{G:U}\abs{G:V}$ which holds when $UV=G.$ (4) is clear (the hypothesis that the extension is central, is only used to guarantee that $G_1$ is $p$-solvable for each prime $p$ dividing $\abs{G_1:\pi^{-1}(U)}=\abs{G:U}$).
\end{proof}

\begin{lem}\label{rcilemma}
Let $G$ be a group and suppose $U=\bigcap_{i=1}^nM_i$ is the regular intersection of maximal subgroups $M_i$ of $G.$ Then the following hold.
\begin{enumerate}
\item $U$ is the complete intersection of the $M_i$ in $G.$
\item $U$ either covers or avoids any chief factor of $G.$
\item Let $M$ be a maximal subgroup of $G$ such that $G$ is $p$-solvable for a prime dividing $\abs{G:M}.$ Then either $MU=G,$ or $M$ is conjugate to exactly one of the $M_i.$
\item The set $\set{M_i}$ consists of exactly one member of each conjugacy class of maximal subgroups of $G$ which contain a
conjugate of $U.$
\end{enumerate}
\end{lem}

\begin{proof}
1. Let $1=K_1\subset K_2\subset...\subset K_m=G$ be any chief series of $G.$ By Lemma \ref{avoidlemma}, each $M_i$ avoids exactly one factor $K_j/K_{j-1},$ and then  $K_j/K_{j-1}\cong X_{M_i}.$ Since no two $X_{M_i}$ are isomorphic, it follows that the different $M_i$ avoid different chief factors. $U$ avoids all these factors, so $\abs{G:U}=\prod_{i=1}^m\abs{K_iU:K_{i-1}U}$ is divisible by $\prod_{i=1}^m\abs{G:M_i}.$ Since it is always the case that $\abs{G:\bigcap_iM_i}\leq\prod_{i=1}^m\abs{G:M_i},$ equality must hold. The relevant solvability conditions are the same for complete and regular intersections, so are satisfied. Hence $U=\bigcap_i M_i$ is a complete intersection, which is (1). 

2. In (1), the remaining factors $\abs{K_iU:K_{i-1}U}$ above must be unity, so $U$ covers all those chief factors which are not avoided and (2) follows. 

3. From (1) and (2), in any chief series, $U$ avoids $n$ chief factors, which are the (distinct) chief factors avoided by $M_1,...,M_n$ respectively, and $U$ covers the remaining factors. By Lemma \ref{XMlem}, $C_G(X_M)/\core_G(M)$ is an abelian chief factor of $G$ avoided by $M.$ Choose a chief series of $G$ which includes the section $\core_G(M)\subset C_G(X_M).$ If $C_G(X_M)/\core_G(M)$ is covered by $U,$ then in our chosen series, $M$ avoids the factor $C_G(X_M)/\core_G(M)$ and $U$ avoids $n$ other distinct factors. Therefore $U\cap M$ avoids at least these $n+1$ factors, and by considering their orders, we easily see that $\abs{G:U\cap M}\geq\abs{G:U}\abs{G:M}$ or $$\frac{\abs{U}\abs{M}}{\abs{U\cap M}}\geq\abs{G},$$ so equality holds and $UM=G,$ which is the first alternative. 

Otherwise, just one of the $M_i$ must avoid $C_G(X_M)/\core_G(M).$ We may choose notation such that this is $M_1.$ If $M$ and $M_1$ are not $G$-conjugate, then by Lemma \ref{lem3}, $\core_G(M)$ and $\core_G(M_1)$ are distinct subgroups of $C_G(X_M).$ Also by Lemma \ref{avoidlemma}, $X_{M_1}$ and $X_M$ are isomorphic. Hence $M_1\cap C_G(X_M)=M_1\cap C_G(X_{M_1})= \core_G(M_1),$ and so $M_1\cap\core_G(M)=\core_G(M_1)\cap\core_G(M).$ We conclude that $M_1,$ and hence $U,$ avoids the factor $\core_G(M)/(\core_G(M_1)\cap\core_G(M)).$ Clearly $M$ covers this factor, and also
$$\frac{\core_G(M)}{\core_G(M_1)\cap\core_G(M)}\cong\frac{\core_G(M)\core_G(M_1)}{\core_G(M_1)}=\frac{C_G(X_M)}{\core_G(M_1)}\cong X_M$$ as $G$-modules. Therefore, in any chief series containing the factor $$\core_G(M)/(\core_G(M_1)\cap\core_G(M)),$$ $U$ avoids this factor, and so avoids no other factor isomorphic to $X_M,$ while $M$ avoids the distinct chief factor $C_G(X_M)/\core_G(M),$ which is isomorphic to $X_M.$ Again, $U\cap M$ avoids $n+1$ chief factors and by considering the index of $U\cap M$ exactly as above, we conclude that $MU=G.$ This completes the proof of (3).

4. Clearly each $M_i$ contains $U.$ Conversely, if $M$ contains a conjugate of $U$ then $G$ is $p$-solvable for any prime dividing $\abs{G:M},$ and also $MU\neq G,$ so by (3), $M$ is conjugate to exactly one of the $M_i$ as required.
\end{proof}

\begin{proof}[Proof of Theorem G]
We have $\bar\chi\chi=(1_U)^G$ where $U$ is a regular intersection. Let $V=\bigcap L_i$ be the complete intersection of the maximal subgroups $L_i.$ If we can show that each $L_i$ contains a conjugate of $U,$ then $V$ contains a conjugate of $U$ by Lemma \ref{cilemma}(1). We may therefore assume that $V=L_1$ is maximal. If $UV=G,$ then $\bigcup_{g\in G}V^g=\bigcup_{u\in U}V^u,$ and since $\bar\chi\chi=(1_U)^G$ vanishes off this set, we have $U\subseteq \bigcup_{u\in U}V^u.$ Therefore, for any $x\in U,$ there exists $u\in U$ such that $x^u\in U\cap V.$ However, no proper subgroup of a finite group can meet each conjugacy class, so $U\cap V=U$ and $U\subseteq V,$ as desired. Therefore $UV\neq G.$ But then by Lemma \ref{rcilemma}(3), $V$ contains a conjugate of $U.$ This completes the proof of the first part. 

For the last statement, if $V$ is maximal and does not contain a conjugate of $U,$ then $UV=G$ by Lemma \ref{rcilemma}(3). Then $[\chi,\chi]_V=[(1_U)^G,1]_V=1,$ so $\chi\vert_V$ is irreducible as required.
\end{proof}

\section{Proof of Theorems C and D}
Suppose $G$ is $p$-solvable and $\chi\in\Irr(G)$ is faithful and strongly irreducible, with $p$ dividing $\chi(1).$ Let $Z$ be the
centre of $G$ and let $K/Z$ with $K\supset Z$ be a chief factor of $G.$ Then by Proposition \ref{prop2}(3), $K/Z$ is the unique
minimal normal subgroup of $G/Z,$ and by Corollary \ref{maxcor1}, $G$ has a maximal subgroup $U,$ unique up to conjugacy, satisfying
$G=KU\text{ and }K\cap U=Z.$ In this situation, provided $p$ is odd, we show that $\bar\chi\chi=(1_U)^G.$

\begin{prop}\label{prop1UG}
Let $G,p,\chi,K,Z$ and $U$ be chosen as above, and suppose $p$ is odd. Then $\bar\chi\chi=(1_U)^G.$
\end{prop}

\begin{proof}
First, suppose $u\in U$ and $\chi(u)\neq 0.$ Then by Proposition \ref{prop2}(2), $\bar\chi(u)\chi(u)=\abs{C_{G/Z}(uZ)}.$ On the other hand, $$(1_U)^G(u)=\sharp\set{gU, u\in U\cap U^g}=\frac{1}{\abs{Z}}\sharp\set{k\in K, u\in U\cap U^k}.$$ 
If $u\in U\cap U^k,$ then $u^{-1}kuk^{-1}\in U\cap K=Z$ so $kZ\in C_{G/Z}(uZ),$ and conversely if $kZ\in C_{G/Z}(uZ),$ then $u^k\in uZ,$ so $u\in U\cap U^k.$ Therefore $$(1_U)^G(u)=\abs{C_K(u)/Z}=\bar\chi(u)\chi(u).$$
Now suppose we can show that $u\in U$ implies $\chi(u)\neq 0.$ Then the above reasoning applies to any $u\in U,$ so
$\bar\chi\chi$ and $(1_U)^G$ agree on $U.$ However $(1_U)^G$ is zero off the union of the conjugates of $U,$ and since
$\sum_{g\in G}(1_U)^G(g)=\abs{G}=\sum_{g\in G}\bar\chi(g)\chi(g),$
and $\bar\chi(g)\chi(g)\geq 0$ for all $g\in G,$ so in fact $\bar\chi\chi$ must also vanish off the union of the conjugates of $U,$ and we have $\bar\chi\chi=(1_U)^G.$

It therefore remains to show that $\chi$ does not vanish on $U.$ Choose a transversal $T$ to $Z$ in $G.$ For $g\in G,$ let $g^*\in Z$
be the unique element with $g=tg^*$ and $t\in T.$ Let $f:G/Z \times G/Z\rightarrow Z$ be defined by
$f(g_1Z,g_2Z)=(g_1g_2)^*/(g_1^*g_2^*).$ Then $f$ is well defined and is a $2$-cocycle on $G$ ($f$ is just the factor set
associated with the extension $G$ of $Z$). We have $\chi\vert_Z=\chi(1)\lambda$ where $\lambda:Z\rightarrow \mathbb{C}^\times$ is a linear character. Then $f_\lambda$ defined by
$f_\lambda=\lambda\circ f$ is a $2$-cocycle $f_\lambda:G/Z \times G/Z\rightarrow \mathbb C^\times.$ Hence $f_\lambda$ defines a
class $[f_\lambda]\in H^2(G/Z,\mathbb{C}^\times).$ If $[g_1,g_2]=z\in Z$ then $g_1g_2=g_2g_1z,$ so
$f_\lambda(g_1Z,g_2Z)/f_\lambda(g_2Z,g_1Z)=(g_1g_2)^*/(g_2g_1)^*=\lambda(z).$ Since $2$-coboundaries are symmetric, the above
equation remains true if $f_\lambda$ is replaced by a cohomologous $2$-cocycle. By Corollary \ref{sigmacor} applied to $G/Z,$
$f_\lambda^\sigma$ is cohomologous to $f_\lambda,$ where $\sigma$ is the automorphism of $G/Z$ which acts trivially on $U/Z$ and
inverts $K/Z.$ Hence if $k\in K, u\in U$ and $[k,u]=z\in Z$ then $[k^{-1},u]=z^{-1}$ and
$$\lambda(z)=f_\lambda(kZ,uZ)/f_\lambda(uZ,kZ)=f_\lambda(k^{-1}Z,uZ)/f_\lambda(uZ,k^{-1}Z)=\lambda(z^{-1}).$$
Since $\lambda$ is faithful, we have $z^2=1.$ But since $z=[k,u]$ is central in $G,$ we have $z^p=[k,u]^p=[k^p,u]=1,$ where the last equality is because $k^p\in Z.$ Since $p$ is odd, we have $z=1.$ Hence $C_{K/Z}(u)=C_K(u)/Z.$ By Proposition \ref{prop2}, $\chi$ does
not vanish on $U,$ as required.
\end{proof}

\begin{cor}\label{cor1}
Let $\chi\in\Irr(G)$ and suppose $\chi(1)$ is divisible by the odd prime $p$ and $G$ is $p$-solvable. Then $\chi$ is strongly
irreducible if and only if $\bar\chi\chi=(1_U)^G$ for a maximal subgroup $U\subseteq G.$
\end{cor}

\begin{proof}
The ``only if'' part follows from Proposition \ref{prop1UG} applied to $G/\ker\chi,$ using Lemma \ref{triviallemma}. For the other implication, suppose
$\bar\chi\chi=(1_U)^G$ for a maximal subgroup $U.$ Then $Z(\chi)=\core_G(U),$ so if $N\vartriangleleft G$ then either $NU=G,$ in which
case $\chi\vert_N$ is irreducible, or $N\subseteq\core_G(U)=Z(\chi).$
\end{proof}

We have given an essentially self-contained derivation of Corollary \ref{cor1}, based
on the elementary Proposition \ref{prop2}. By adopting this approach, we have given self-contained proofs of Theorems A and B along
the way. The most difficult aspect of Corollary \ref{cor1} is the fact that $\chi$ does not vanish on $U,$ and in our treatment this
follows from properties of the Schur multiplier of $G/Z$ given by Proposition \ref{MGprop} and Corollary \ref{sigmacor}.

However, Corollary \ref{cor1} can alternatively be inferred from a fundamental theorem of Isaacs on
ramified sections. Before proceeding, we show how this is done. Recall that a \emph{character five} $(G,K,L,\theta,\phi)$ is a group
$G$ with normal subgroups $K,L\vartriangleleft G$ with $L\subseteq K,$ $\theta\in\Irr(K),$ $\phi\in\Irr(L)$ such that $\phi$ is a
constituent of $\theta_L$ and $\theta,\phi$ are invariant in $G$ and fully ramified with respect to $K/L.$

\begin{thm}[\cite{Isaacs1973}, Proposition 9.1]\label{bigthm}
Let $(G,K,L,\theta,\phi)$ be a character five and assume that either $\abs{G:K}$ or $\abs{K/L}$ is odd. Then there is a (reducible)
character $\psi\in\mathrm{Char}(G)$ and a subgroup $U\subset G$ with the properties that
\begin{enumerate}
\item $UK=G$ and $U\cap K=L.$
\item For $g\in G,$ $\abs{\psi(g)}^2=\abs{C_{K/L}(g)}$
\item If $g\in G$ is not $G$-conjugate to an element of $U$ then $\chi(g)=0$ for all $\chi\in\Irr(G\vert\theta).$
\item The equation $\chi\vert_U=\psi\vert_U\xi,$ for $\chi\in\Irr(G\vert\theta)$ and $\xi\in\Irr(U\vert\phi)$ defines a one-to-one
correspondence between these sets of characters.
\item If $\abs{G:L}$ is odd, then $\chi$ and $\xi$ correspond above, if and only if $2\nmid [\chi\vert_U,\xi].$
\end{enumerate}
\end{thm}

\begin{cor}\label{bigcor}
In the situation of Theorem \ref{bigthm}, let $\chi\in\Irr(G\vert\theta)$ and suppose that $\phi(1)=1$ and $\chi\vert_K$ is
irreducible. Then
\begin{enumerate}
\item $\bar\chi\chi=(1_U)^G.$
\item $\chi\vert_U$ has a linear constituent.
\end{enumerate}
\end{cor}

\begin{proof}
By hypothesis, $\chi\vert_K$ is irreducible, so since $\chi\in\Irr(G\vert\theta),$ we have $\chi\vert_K=\theta.$ By the hypotheses of
Theorem \ref{bigthm}, $(\theta(1)/\phi(1))^2=\abs{K/L},$ so in fact $\chi(1)^2=\abs{K/L}.$ However by Theorem \ref{bigthm}(2),
$\psi(1)^2=\abs{K/L},$ and we conclude that $\chi(1)=\psi(1).$ But by Theorem \ref{bigthm}(3), $\chi\vert_U=\psi\vert_U\xi,$ where
$\xi\in\Irr(U).$ Hence $\xi$ is linear. Therefore by Theorem \ref{bigthm}(2), (3),

$$\abs{\chi(g)}^2=\abs{\psi(g)}^2=\begin{cases}\abs{C_{K/L}(g)} & \text{ if } g\in U,\\$0$ & \text{ otherwise.}\end{cases}$$

As in the first part of the proof of Proposition \ref{prop1UG}, the right hand side of the above equation is just $(1_U)^G(g),$ which
proves (1). Also, by Theorem \ref{bigthm}(5), $\chi\vert_U$ contains the linear constituent $\xi,$ which gives (2).
\end{proof}

\begin{proof}[Proof of Theorem C]
First, we show that we can reduce to the case $\chi$ strongly irreducible. By Proposition \ref{prop1}, there exist a central extension
$\pi:\tilde G\twoheadrightarrow G$ and strongly irreducible characters $\rho_1,...,\rho_n\in\Irr(\tilde G)$ with
$\prod_{i=1}^n\rho_i=\tilde\chi,$ where $\tilde\chi=\chi\circ\pi$ is the lift of $\chi$ to $\tilde G.$ By Corollary \ref{cor1}, for
each $i$ there is a maximal subgroup $\tilde M_i\subseteq G$ with $\bar\rho_i\rho_i=(1_{\tilde M_i})^{\tilde G}.$ Assuming that Theorem C holds for strongly irreducible characters, we have that $\mathfrak M_0(\rho_i)=\tilde M_i^{\tilde G}.$ For each $i,$
$\core_{\tilde G}\tilde M_i= Z(\rho_i)\supseteq Z(\tilde G)\supseteq\ker\pi,$ so $M_i=\pi(\tilde M_i)$ is a maximal subgroup of $G,$ and clearly
$\bar\rho_i\rho_i$ may be viewed as a character of $G$ and then has $\bar\rho_i\rho_i=(1_{M_i})^G.$ Put $U=\bigcap_i M_i.$ Then $$[\prod_i(1_{M_i})^G,1]_G=[\bar\chi\chi,1]_G=1,$$ so $G$ is transitive on $\prod_i(1_{M_i})^G.$ It
follows that $U$ is the complete intersection of the $M_i.$ Also $\bar\chi\chi=(1_U)^G.$ If $V\subseteq G$ is a complete
intersection and $\chi\vert_V$ has no zeros, then by Lemma \ref{cilemma}(4), $\pi^{-1}(V)$ is a complete intersection in $\tilde G.$ No $\rho_i$ has a zero on $\pi^{-1}(V),$ so since we are assuming that Theorem C holds for each $\rho_i,$ $\pi^{-1}(V)$ is conjugate
to a subgroup of $\tilde M_i$ for each $i.$ By Lemma \ref{cilemma}(1), $V$ is conjugate to a subgroup of $U,$ as required.

We are left with the case that $\chi$ is strongly irreducible. We may assume $\chi$ is non-linear, since otherwise the result is trivial with $\mathfrak M_0(\chi)=\set{G}.$ By Corollary \ref{cor1} again, there is a maximal subgroup $U\subseteq G$ with
$\bar\chi\chi=(1_U)^G.$ Clearly $U\in\mathfrak M_0(\chi).$ Now suppose $L=\bigcap_{i=1}^nM_i\subseteq G$ is the complete intersection of the maximal subgroups $M_i\subset G,$ and suppose $\chi\vert_L$ has no zeros. We have to show that $L\subseteq_G U.$ We use the partial ordering introduced below Lemma \ref{lem3}, noting that since $L=\bigcap_{i=1}^nM_i\subseteq G$ is a complete intersection, $G$ is $p$-solvable for each prime $p$ dividing $\abs{G:M_i}.$ First, suppose there exists $i$ such that $M_i^G\ngeq M_j^G$ for all $j\neq i,$ and $M_i^G\ngeq U^G.$ We may choose notation such that $i=1.$ Then $L=\bigcap_{i=2}^n(M_1\cap M_i),$ and by Corollary \ref{maxsubgroupcor}, each $M_1\cap M_i$ is a maximal subgroup of
$M_1.$ Also since $L$ is a complete intersection in $G,$ by Lemma \ref{cilemma}(2), $M_1M_i=G$ for all $i\neq 1,$ and so from $$\abs{G:L}=\prod_{i=1}^n\abs{G:M_i},$$ we have $$\abs{M_1:L}=\prod_{i=2}^n\abs{G:M_i}=\prod_{i=2}^n\abs{M_1M_i:M_i}=\prod_{i=2}^n\abs{M_1:M_1\cap M_i},$$ which shows that $L$ is a complete intersection in $M_1.$ Also, since $M_1^G\ngeq U^G,$ $U\cap M_1$ is a maximal subgroup of $M_1$ and $UM_1=G,$ so $\bar\chi\chi\vert_{M_1}=(1_{U\cap M_1})^{M_1}$ and by Corollary \ref{cor1}, $\chi\vert_{M_1}$ is strongly irreducible. Hence $U\cap M_1\in\mathfrak M_0(\chi\vert_{M_1}).$ By induction, $L$ is $M_1$-conjugate to a subgroup of $U\cap M_1,$ so $L$ is certainly $G$-conjugate to a subgroup of $U,$ as required. Hence we may assume that no such $i$ exists, so for every $i,$ at least one of the following holds.
\begin{enumerate}
\item There exists $j\neq i$ such that $M_i^G\geq M_j^G.$
\item $M_i^G\geq U^G.$
\end{enumerate}
Since $L$ is a complete intersection, $M_iM_j=G$ for $i\neq j$ by Lemma \ref{cilemma}(2), so $M_i^G\neq M_j^G.$ Therefore, if (1) holds for any $i,$ then $M_i^G>M_j^G,$ and repeating with $M_j,$ eventually (2) must hold. Hence in fact (2) holds for each $i,$ and we conclude that $Z(\chi)=\core_G(U)\subseteq\core_G (M_i)$ for each $i.$ If no $M_i$ is conjugate to $U,$ then by Lemmas \ref{lem3} and \ref{XMlem}, we have $C_G(X_U)\subseteq L.$ But since $\chi$ is strongly irreducible, $\chi\vert_{C_G(X_U)}$ is irreducible. Then since $\chi$ is non-linear, $\chi$ must have a zero on $C_G(X_U).$ This conflicts with $C_G(X_U)\subseteq L,$ and this contradiction completes the proof.
\end{proof}

The proof of Theorem D is free, being essentially contained in the first part of the proof of Theorem C.

\begin{proof}[Proof of Theorem D]
The first part of the proof of Theorem C gave a complete intersection $U=\bigcap_i M_i$ with $\bar\chi\chi=(1_U)^G$ and the statements of Theorem D all satisfied. By Theorem C, all elements of $\mathfrak M_0(\chi)$ are conjugate and so enjoy the same properties.
\end{proof}

\section{Proof of Theorems E, F, J and K}

\begin{lem}\label{goinguplemma}
Let $K/L$ be an abelian chief factor of $G$ and suppose $\varphi\in\Irr(L)$ is invariant in $G.$ Then either \begin{enumerate}
\item $\varphi$ extends to $K,$ or
\item $\varphi$ is fully ramified with respect to $K/L.$
\end{enumerate}
\end{lem}

\begin{proof}
This is to be found at \cite[Lemma 3.4]{Isaacs1981}, or as part of \cite[Problem 6.12]{Isaacs1976}.
\end{proof}

The following ``going up'' theorem is the key ingredient of Theorem E. Results of this type can be proved using Glauberman's lemma under coprimeness assumptions (see particularly Theorem 13.27 in \cite{Isaacs1976}). This approach does not seem to be sufficient for our application here, but we can succeed by using a short cohomological argument instead.

\begin{thm}\label{goingupthm}
Suppose $G$ is $p$-solvable for a prime $p,$ and let $K/L$ be a $p$-chief factor of $G.$ Suppose $\varphi\in\Irr(L)$ is $G$-invariant, but that no character in $\Irr(K\vert\varphi)$ is $G$-invariant. Let $\eta\in\Irr(K\vert\varphi).$ Then $Irr(K\vert\varphi)$ consists of the conjugates of
$\eta,$ and the inertia group $I_G(\eta)$ is a maximal subgroup of $G.$ Furthermore, $X_{I_G(\eta)}$ is defined and $X_{I_G(\eta)}\cong (K/L)^*$ as $G$-modules.
\end{thm}

\begin{proof}
If $\varphi$ were fully ramified with respect to $K/L,$ then we would have $\varphi^K=e\eta$ with $e^2=\abs{K:L},$ so that $\eta$ would be $G$-invariant, contrary to hypothesis. Hence by Lemma \ref{goinguplemma}, we may assume $\eta$ is an extension of $\varphi.$ Therefore, by \cite[Corollary 6.17]{Isaacs1976}, the characters $\eta\lambda$ are distinct for $\lambda\in\Irr(K/L),$ and are all the irreducible characters of $K$ lying over $\varphi.$ If $g\in G$ then $\eta^g$ also lies over $\varphi,$ so $\eta^g=\eta\lambda_g$ for a unique $\lambda_g\in\Irr(K/L).$ Now if $g,h\in G$ then $$\eta^{gh}=\eta\lambda_{gh}=(\eta\lambda_g)^h=\eta\lambda_h\lambda_g^h,$$ so $\lambda_{gh}=\lambda_h\lambda_g^h.$ Writing $X$ for the $G$-module $K/L,$ let $\delta:G\rightarrow X^*$ be the map which takes $g\mapsto\lambda_g.$ The preceding sentence says that $\delta$ is a derivation from $G$ into the right $G$-module $X^*.$ Suppose $\delta$ were inner. Then there would be $\mu\in\Irr(X)=X^*$ with $\lambda_g=\mu^g\mu^{-1}.$ Then $(\eta\mu^{-1})^g=\eta\lambda_g\mu^{-g}=\eta\mu^{-1},$ which shows that $\eta\mu^{-1}\in\Irr(K\vert\varphi)$ is $G$-invariant, contrary to hypothesis. Therefore $[\delta]\in H^1(G,X^*)$ cannot be zero, so since $G$ is $p$-solvable and $C_G(X)=C_G(X^*),$ Lemma \ref{firstcohomologylemma} shows that $$\epsilon=\res_{C_G(X)}(\delta)\in\hm_G(C_G(X),X^*)\neq 0.$$ Since $\epsilon$ is not zero and $X$ is simple, it is surjective. Hence if $\lambda\in\Irr(K/L)$ then $\lambda=\lambda_h$ for some $h\in C_G(X),$ and then $\eta^h=\eta\lambda.$ In particular, every character in $\Irr(K\vert\varphi)$ is conjugate to $\eta$ by an element of $C_G(X).$ Hence they are certainly all conjugate in $G,$ which proves the first statement, and moreover $C_G(X)I_G(\eta)=G.$ Now $$\ker\epsilon=\set{g\in C_G(X)\vert\lambda_g=1}=C_G(X)\cap I_G(\eta),$$ and also $$\core_G(I_G(\eta))=\bigcap_{g\in G} I_G(\eta\lambda_g)= C_G(X^*)\cap I_G(\eta)=C_G(X)\cap I_G(\eta).$$ Hence $\ker\epsilon=\core_G(I_G(\eta)),$ so $C_G(X)/\core_G(I_G(\eta))$ is an abelian chief factor of $G$ isomorphic to $X^*.$ Then in $G/\core_G(I_G(\eta)),$ $I_G(\eta)$ is a complement to the abelian minimal normal subgroup $C_G(X)/\core_G(I_G(\eta)).$ It follows that $I_G(\eta)$ is maximal and has $X_{I_G(\eta)}\cong (K/L)^*,$ as required.
\end{proof}

\begin{proof}[Proof of Theorem E]
We have $\chi\in\Irr(G)$ with $\bar\chi\chi=(1_U)^G,$ where $U\subseteq G$ is a regular intersection, and we must show that $\chi$ is quasi-primitive. Suppose not, and let $N\vartriangleleft G$ be such that $\chi\vert_N$ is not homogeneous. Choose a chief series through $N,$ and let $L\subset N$ be the largest member of the series for which $\chi\vert_L$ is homogeneous. Let $K$ be the next member up, so that $K/L$ is a chief factor of $G,$ $\chi\vert_L$ is homogeneous, but $\chi\vert_K$ is not homogeneous. Let $\eta\in\Irr(K)$ and $\varphi\in\Irr(L)$ be constituents of $\chi\vert_K$ and $\chi\vert_L$ respectively. We have $$\chi\vert_K=e\sum_{g\in G/I_G(\eta)}\eta^g\text{ and }\chi\vert_L=f\varphi,$$ for suitable multiplicities $e,f.$ 

We claim that $U$ avoids $K/L.$ First, note that $$[\chi,\chi]_K=[(1_U)^G,1]_K=\abs{G:UK}$$ and similarly $[\chi,\chi]_L=\abs{G:UL}.$
If $U$ does not avoid $K/L$ then by Lemma \ref{rcilemma}(2), $U$ covers $K/L$ and we have $[\chi,\chi]_K=[\chi,\chi]_L.$ Now $$e^2\abs{G:I_G(\eta)}=[\chi,\chi]_K=f^2\text{ and }e\abs{G:I_G(\eta)}\eta(1)=\chi(1)=f\varphi(1),$$ and dividing the second of these two equations by the first on both sides gives $\eta(1)/e=\varphi(1)/f,$ or $\eta(1)/\varphi(1)=e/f.$ However, $\varphi(1)$ divides $\eta(1)$ and $e$ divides $f.$ Hence $e=f$ and $\eta(1)=\varphi(1),$ and we conclude from the first relation that $\abs{G:I_G(\eta)}=1.$ But this implies that $\chi\vert_K$ is homogeneous, which is not the case. Therefore $U$ avoids $K/L$ as claimed. In particular, by Lemma \ref{calem}, $\abs{K/L}$ divides $\abs{G:U}.$ Since by definition of a regular intersection, $G$ is $\pi(\abs{G:U})$-solvable, therefore $K/L$ is a $p$-chief factor for some prime $p$ and $G$ is $p$-solvable. 

Also, since $\bar\chi\chi\vert_K=(1_U)^G\vert_K$ and $U$ avoids $K/L,$ we find that $\chi$ vanishes on $K\backslash L.$ Therefore if $\theta\in\Irr(K\vert\varphi),$ then $[\chi,\theta]_K\neq 0,$ so $\theta=\eta^g$ for some $g\in G.$ In particular, $\theta$ cannot
be $G$-invariant. Now Theorem \ref{goingupthm} applies, and $I_G(\eta)$ is a maximal subgroup of $G$ with $X_{I_G(\eta)}\cong (K/L)^*.$ Also, since $K\subseteq\core_G(I_G(\eta)),$ we may choose a chief series passing through $L,K$ and $\core_G(I_G(\eta)).$ Since both $I_G(\theta)$ and $U$ are regular intersections, by Lemma \ref{rcilemma}(2), they cover or avoid each factor, and in the chosen series, $U$ and $I_G(\eta)$ avoid distinct
factors. By considering the orders of $U,I_G(\eta)$ and $U\cap I_G(\eta),$ it follows that $UI_G(\eta)=G.$ Hence $[\chi,\chi]_{I_G(\eta)}=[(1_U)^G,1]_{I_G(\eta)}=1,$ so $\chi\vert_{I_G(\eta)}\in\Irr({I_G(\eta)})$ and $\chi\vert_K$ is homogeneous. This
contradiction completes the proof.
\end{proof}

We next give the proof of Theorem F. This is based on similar ideas to the proof of \cite[Proposition 1.4]{FergusonTurull1985} in the odd order case, and makes the same use of a result of E. C. Dade on anisotropic modules for groups of odd order (\cite[Corollary 2.10 and Proposition 1.10]{Dade1981}).

\begin{proof}[Proof of Theorem F]
Let $U\in\mathfrak M_0(\chi),$ and suppose $M_1,M_2$ are non-conjugate maximal subgroups such that $X_{M_1}$ is isomorphic to at least one of $X_{M_2}$ or $X^*_{M_2}$ as $G$-modules, and $U\subseteq M_1\cap M_2.$ In particular, $G$ is $p$-solvable for each prime dividing $\abs{G:M_1}=\abs{G:M_2}.$ Let $K=C_G(X)$ and $L=\core_G(M_1)\cap\core_G(M_2).$ Since $M_1$ and $M_2$
are not conjugate, $\core_G(M_1)$ and $\core_G(M_2)$ are distinct by Lemma \ref{lem3}. Since $K/\core_G(M_1)$ is a chief factor of $G$ and $\core_G(M_1)\subset\core_G(M_1)\core_G(M_2)\subseteq K,$ we must have $\core_G(M_1)\core_G(M_2)=K,$ and therefore $$K/L=K/\core_G(M_1)\times K/\core_G(M_2).$$ Hence $K/L$ is elementary abelian, and as $G$-module (written multiplicatively) we have $K/L\cong X_{M_1}\times X_{M_2}.$

Since $\chi$ is quasi-primitive, $\chi\vert_L=f\theta,$ for some $G$-invariant $\theta\in\Irr(L).$ Suppose $L\subseteq N\subseteq K$ with
$N\vartriangleleft G,$ then similarly $\chi\vert_N=e\varphi$ for some $\varphi\in\Irr(N)$ lying over $\theta.$ Now $\chi$ vanishes on
$N\backslash L,$ hence so does $\varphi.$ Hence $\theta$ is fully ramified with respect to $N/L.$

We consider the form $\ll x,y\gg_\theta$ defined for the character triple $(K,L,\theta)$ (see \cite[Problem 11.12]{Isaacs1976} or
\cite[Section 2]{Isaacs1973}). By \cite[Lemma 2.7]{Isaacs1973}, since $\theta$ is fully ramified with respect to any
$N\vartriangleleft G$ with $L\subseteq N\subseteq K,$ then $\ll x,y\gg_\theta$ is anisotropic (that is, $K/L,$ viewed as $G$-module, has no non-trivial isotropic $G$-submodule). The restriction of $\ll x,y\gg_\theta$ to any $G$-submodule of $K/L$ is then also anisotropic, and therefore non-singular. Hence in fact $X_1\cong X^*_1,$ and it follows that $X_1$ and $X_2$ are isomorphic. Hence $K/L$ is a homogeneous, not simple, anisotropic $G$-module. By \cite[Corollary 2.10 and Proposition 1.10]{Dade1981}, it follows that $\abs{G:K}$ is even. Since $Z(\chi)=\core_G(U)\subseteq K,$ we then have $\abs{G:Z(\chi)}$ even. This completes the proof.
\end{proof}

For Theorem J, we need the additional information provided by the following theorem of Cossey, Hawkes and Willems.

\begin{thm}\label{CHWthm}
Let $G$ be $p$-solvable and let $H$ be a Hall $p$-complement in $G.$ Suppose $X$ and $Y$ are non-isomorphic simple $\mathbb F_pG-modules,$ whose dimensions over $\mathbb F_p$ are prime to $p.$ Then $X\vert_H$ and $Y\vert_H$ are non-isomorphic simple $\mathbb F_pH-modules.$
\end{thm}

\begin{proof}
This is immediate from Theorems 1 and 2 of \cite{CosseyHawkesWillems1980}, noting that the dimension of any absolutely irreducible constituent of $X$ divides the dimension of $X,$ and so is prime to $p,$ and similarly for $Y.$ 
\end{proof}

\begin{proof}[Proof of Theorem J]
We have $\bar\chi\chi=(1_U)^G$ with $U$ the regular intersection of the maximal subgroups $L_i,$ say. Note that $G$ is $p$-solvable because $p$ divides $\abs{G:U}.$ Since $U$ is a complete intersection by Lemma \ref{rcilemma}(1), we have $$\chi(1)^2=\abs{G:U}=\prod_i\abs{X_{L_i}}=p^{\sum_i\dim(X_{L_i})}=p^{2a}.$$ Therefore $\sum_i\dim(X_{L_i})<2p.$ By Theorem D(2), moreover, each $\abs{X_{L_i}}$ is a square, or equivalently, each $\dim(X_{L_i})$ is even. In particular, $\dim(X_{L_i})\neq p,$ so $\dim(X_{L_i})$ is prime to $p$ for each $i.$ Since $H$ contains a Hall $p$-complement of $G,$ by Theorem \ref{CHWthm}, the $X_{L_i}$ remain irreducible on restriction to $H,$ and moreover if either $X_{L_i}\cong X_{L_j}$ or $X_{L_i}\cong X^*_{L_j}$ as $H$-modules, then the same holds as $G$-modules, and so $i=j.$

Since $\chi\vert_H$ is irreducible, we have $[(1_U)^G,1]_H=1$ and it follows that $HU=G.$ Now by Proposition \ref{goingdownprop}, each $H\cap L_i$ is a maximal subgroup of $H$ with $X_{H\cap L_i}\cong X_{L_i}\vert_H$ as $H$-modules. Therefore $H\cap U$ is a regular intersection in $H.$ But $\bar\chi\chi\vert_H=(1_U)^G\vert_H=(1_{H\cap U})^U,$ so by Theorem E, $\chi\vert_H$ is quasi-primitive, as required.
\end{proof}

\begin{proof}[Proof of Theorem K]
We have $\chi\in\Irr(G),$ $G$ is $\pi(\chi(1))$-solvable and $\abs{G:Z(\chi)}$ is odd. Suppose $H\subseteq G$ and $\theta$ is a
character of $H$ with $\theta^G=e\chi$ for some $e\in\mathbb N.$ Suppose for a contradiction that this is possible with $H\neq G.$ 

By Theorem D, there is a central extension $\pi:\tilde G\rightarrow
G$ such that the lift $\tilde\chi=\chi\circ\pi$ factors into a product of strongly irreducible characters of $\tilde G.$ Let $\tilde
H=\pi^{-1}(H),$ the pre-image of $H$ in $\tilde G$ and $\tilde\theta=\theta\circ\pi,$ the lift of $\theta$ to $\tilde H.$ Then
$\tilde\theta^{\tilde G}=e\tilde\chi$ and $\tilde H\neq\tilde G,$ so we may replace $G,H,\chi,\theta$ with $\tilde G,\tilde
H,\tilde\chi,\tilde\theta$ and hence assume that $\chi=\prod_{i=1}^n\rho_i$ is a product of strongly irreducible characters in $G.$ We
may also clearly assume that $H$ is maximal and $\theta$ is irreducible.

Let $M_i\in\mathfrak M_0(\rho_i).$ Then $U=\bigcap_iM_i\in\mathfrak M_0(\chi)$ and by Theorem F, $U$ is the regular intersection of
the $M_i.$ Since $e\chi=\theta^G,$ $\chi$ vanishes off conjugates of $H,$ so by Lemma \ref{rcilemma}(3), $H$ is conjugate to exactly
one of the $M_i.$ By choosing notation, and replacing $H$ with a conjugate subgroup if necessary, we may assume that $H=M_1.$ Let
$\psi=\prod_{i=2}^n\rho_i,$ so that $\chi=\rho_1\psi,$ and let $V=\bigcap_{i=2}^nM_i.$ Then $\bar\psi\psi=(1_V)^G$ and since $\chi$ is
irreducible, we have $[(1_H)^G,(1_V)^G]=1,$ so that $HV=1$ and $\psi\vert_H$ is irreducible. Now from
$e\chi(1)=e\rho_1(1)\psi(1)=\abs{G:H}\theta(1)$ and $\rho_1(1)^2=\abs{G:H},$ we deduce that $$e\psi(1)=\theta(1)\rho_1(1),$$ and by
Frobenius reciprocity,$$[\theta,\chi]_H=e=[\theta,\rho_1\psi]_H=[\theta\bar\rho_1,\psi]_H.$$ Since $\psi\vert_H$ is irreducible, we must
have $$\theta\bar\rho_1\vert_H=e\psi\vert_H.$$
At this point, we can observe that $e>1,$ giving the odd order case of Berger's theorem. This follows because $\rho_1\vert_H$ is not irreducible (to see this, for example, we can use $\bar\rho_1\rho_1=(1_H)^G$ to check that $[\rho_1,\rho_1]_H=[\bar\rho_1\rho_1,1]_H>1$). However, we do not use this. Instead, Corollary \ref{bigcor}(2) shows that $\rho_1\vert_H$ has a linear constituent $\xi.$ Therefore,
since $\theta$ and $\psi\vert_H$ are irreducible, $$\theta\bar\xi=\psi\vert_H.$$ Therefore $\theta(1)=\psi(1)$ and
$$e^2=\abs{G:H}^2\psi(1)^2/\chi(1)^2=\abs{G:H}^2/\rho_1(1)^2=\abs{G:H}=\abs{HV:H}=\abs{V:U}.$$ Also, $e\chi=\theta^G=(\psi\vert_H\xi)^G=\psi\xi^G,$ and, using $\bar\psi\psi=(1_V)^G,$ we have
$$e^2=[\theta^G,\theta^G]_G=[\psi\xi^G,\psi\xi^G]_G=[\bar\psi\psi,\bar\xi^G\xi^G]_G=[\xi^G\vert_V,\xi^G\vert_V]_V=[(\xi\vert_U)^V,(\xi\vert_U)^V]_V,$$
so $$e^2=\sum_{g\in U\backslash V/U}[\xi^g,\xi]_{U^g\cap U}.$$ Since $\xi$ is linear, we have $[\xi^g,\xi]_{U^g\cap U}\leq 1,$ so
$e^2=\abs{V:U}\leq\abs{U\backslash V/U},$ so equality must hold, and we conclude that $U\vartriangleleft V$ and $\xi\vert_U$ is
invariant in $V.$ Now since $HV=G$ and $\bar\rho_1\rho_1=(1_H)^G,$ we have $[\rho_1\vert_V,\rho_1\vert_V]=[(1_H)^G,1]_V=1$ and
$\rho_1\vert_V$ is irreducible. Therefore since $U\vartriangleleft V,$ and $\rho_1\vert_U$ contains the $V$-invariant character
$\xi_U,$ so $\rho_1\vert_U=\rho_1(1)\xi\vert_U$ and we conclude that $U\subseteq Z(\rho_1)=\core_G(H).$ Therefore as
$\bar\chi\chi=(1_U)^G,$ $\chi$ vanishes off the normal subgroup $K=\core_G(H).$ Let $N$ be a maximal normal subgroup of $G$ containing
$K.$ Then $G/N$ is cyclic of prime order and $\chi$ vanishes off $N,$ so $\chi\vert_N$ cannot be irreducible, and the only possibility is
that $\chi\vert_N$ is a sum of distinct characters, contradicting the fact that $\chi$ is quasi-primitive. This completes the proof.
\end{proof}

\emph{London, England. Email: tom@beech84.fsnet.co.uk.}
\end{document}